\pdfoutput=1
\RequirePackage{ifpdf}
\ifpdf %We are running pdfTeX in pdf mode
\documentclass[pdftex]{sigma}
\else
\documentclass{sigma}
\fi

\newcommand{\hyp}[5]{{}_{#1}F_{#2} \left(\begin{array}{@{}c@{}} {#3} \\ {#4} \end{array}\!;{#5}\right)}

\numberwithin{equation}{section}
\numberwithin{theorem}{section}
\numberwithin{corollary}{section}
\numberwithin{example}{section}
\numberwithin{note}{section}

\DeclareMathOperator{\vol}{vol}
\DeclareMathOperator{\sn}{sn}

\begin{document}

\newcommand{\arXivNumber}{1405.4847}

\allowdisplaybreaks

\renewcommand{\thefootnote}{$\star$}

\renewcommand{\PaperNumber}{015}

\FirstPageHeading

\ShortArticleName{Fourier and Gegenbauer Expansions on the Hypersphere}

\ArticleName{Fourier and Gegenbauer Expansions\\ for a~Fundamental Solution of Laplace's Equation\\ in Hyperspherical Geometry\footnote{This
paper is a~contribution to the Special Issue on Exact Solvability and Symmetry Avatars in honour of Luc Vinet.
The full collection is available at
\href{http://www.emis.de/journals/SIGMA/ESSA2014.html}{http://www.emis.de/journals/SIGMA/ESSA2014.html}}}

\Author{Howard S.~COHL~$^\dag$ and Rebekah M.~PALMER~$^\ddag$}

\AuthorNameForHeading{H.S.~Cohl and R.M.~Palmer}

\Address{$^\dag$~Applied and Computational Mathematics Division, National Institute of Standards
\\
\hphantom{$^\dag$}~and Technology, Gaithersburg, MD, 20899-8910, USA} \EmailD{\href{mailto:howard.cohl@nist.gov}{howard.cohl@nist.gov}}
\URLaddressD{http://www.nist.gov/itl/math/msg/howard-s-cohl.cfm}

\Address{$^\ddag$~Department of Mathematics, Johns Hopkins University, Baltimore, MD 21218, USA}
\EmailD{\href{mailto:rmaepalmer4@gmail.com}{rmaepalmer4@gmail.com}}

\ArticleDates{Received May 20, 2014, in f\/inal form February 09, 2015; Published online February 14, 2015}

\Abstract{For a~fundamental solution of Laplace's equation on the~$R$-radius~$d$-dimensional hypersphere, we compute the azimuthal
Fourier coef\/f\/icients in closed form in two and three dimensions.
We also compute the Gegenbauer polynomial expansion for a~fundamental solution of Laplace's equation in hyperspherical geometry in
geodesic polar coordinates.
From this expansion in three-dimensions, we derive an addition theorem for the azimuthal Fourier coef\/f\/icients of a~fundamental solution
of Laplace's equation on the 3-sphere.
Applications of our expansions are given, namely closed-form solutions to Poisson's equation with uniform density source distributions.
The Newtonian potential is obtained for the 2-disc on the 2-sphere and 3-ball and circular curve segment on the 3-sphere.
Applications are also given to the superintegrable Kepler--Coulomb and isotropic oscillator potentials.}


\Classification{31C12; 32Q10; 33C05; 33C45; 33C55; 35J05; 35A08; 42A16}

\begin{flushright}
\begin{minipage}{60mm}
\it Dedicated to Professor Luc Vinet,\\
and to his work on the $2$-sphere
\end{minipage}
\end{flushright}

\renewcommand{\thefootnote}{\arabic{footnote}}
\setcounter{footnote}{0}

\section{Introduction}

This investigation is concerned with the study of eigenfunction expansions for a~fundamental solution of Laplace's equation on
Riemannian spaces of positive constant curvature.
These manifolds are classif\/ied as compact Riemannian, and we refer to them as~$d$-dimensional~$R$-radius hyperspheres.
Analysis on hyperspheres has a~long history and has been studied by luminaries such as Schr{\"o}dinger
in 1938 and 1940~\cite{Schrodinger38,Schrodinger40}.
However, there is signif\/icant activity in the literature regarding analysis on hyperspheres
(e.g.,~\cite{GenestVinetZhedanov14,KalMilPost13}, see~\cite{MillerPostWin13} for a~topical review in classical and quantum
superintegrability).
In fact,~\cite{KalMilPost13} makes a~connection between the Askey scheme of hypergeometric orthogonal polynomials and superintegrable
systems on the 2-sphere!

In this paper, we derive expansions for a~fundamental solution of Laplace's equation on the~$d$-dimensional~$R$-radius hypersphere which
allow one to ef\/fectively obtain density-potential pairs for density distributions with a~large degree of rotational symmetry.
Through these expansions, one may further explore classical potential theory and related theories on hyperspheres.
We have previously studied Fourier and Gegenbauer polynomial expansions for a~fundamental solution of Laplace's equation on
the~$d$-dimensional~$R$-radius constant negative curvature non-compact Riemannian manifold given by the hyperboloid model of hyperbolic
geometry~\cite{CohlKalII}.
In the hyperboloid model, Fourier expansions for a~fundamental solution of Laplace's equation were given in rotationally-invariant
coordinate systems and Gegenbauer polynomial expansions of this fundamental solution were given in standard hyperspherical coordinates.
In a~previous paper~\cite{Cohlhypersphere}, we have studied various special function representations for a~fundamental solution of
Laplace's equation on the compact positive constant curvature Riemannian manifold given by the~$d$-dimensional~$R$-radius hypersphere.
In this paper, we derive for this manifold, Fourier and Gegenbauer expansions analogous to those presented in~\cite{CohlKalII}.
This exercise is motivated by the necessity and utility of these expansions in obtaining the solution to potential theoretic problems in
this geometry.

As an application and demonstration of the utility of our expansions, in Section~\ref{ApplicationsinPotentialTheory} we give three
examples in Newtonian potential theory on the 2-sphere and 3-sphere which reduce to elementary results in Euclidean potential theory.
The density distributions we treat possess spherical and rotational symmetry, and their potential solutions are obtained through the
derived Gegenbauer and azimuthal Fourier expansions respectively.
These density distributions are the 2-disc on the 2-sphere, the 3-ball on the 3-sphere, and the rotationally-invariant circular curve
segment on the 3-sphere.
Using superintegrable potentials, one may derive a~paired density distribution assuming that the potentials satisfy Poisson's equation.
As such, one may convolve a~fundamental solution of Laplace's equation on the hypersphere with this paired density distribution to
re-obtain the superintegrable potential.
Examples of this connection to superintegrable potentials are given in Section~\ref{Applicationstosuperintegrable} for the isotropic
oscillator potential and for the Kepler--Coulomb potential on~$d$-spheres, for $d=2,3,\ldots$.

We denote the sets of real and complex numbers by ${\mathbf R}$ and ${\mathbf C}$, respectively.
Similarly, the sets ${\mathbf N}:=\{1,2,3,\ldots\}$ and ${\mathbf Z}:=\{0,\pm1,\pm2,\ldots\}$ denote the natural numbers and the
integers.
Furthermore, we denote the set ${\mathbf N}_0:=\{0,1,2,\ldots\}={\mathbf N} \cup \{0\}$.
Throughout this paper, we adopt the common convention that empty sums and products are zero and unity respectively.
Note that we often adopt a~common notation used for fundamental solution expansions, namely if one takes $a,a^\prime\in{\mathbf R}$, then
\begin{gather}
\label{alessgtr}
a_\lessgtr:={\min \atop \max}\{a,a^\prime\}.
\end{gather}

\section{Preliminaries}

\subsection{Coordinates}\label{Coordinates}

The hypersphere $\mathbf{S}_R^d$, a~smooth~$d$-dimensional compact constant positive curvature Riemannian manifold can be def\/ined as
Riemannian submanifold given by set of all points $\mathbf{x}=(x_0,\ldots,x_d)$ such that $x_0^2+\dots+x_d^2=R^2$.
Standard hyperspherical coordinates are an example of geodesic polar coordinates on this manifold given~by
\begin{gather}
\begin{split}
& x_0=R\cos\theta,
\\
& x_1=R\sin\theta\cos\theta_{d-1},
\\
& x_2=R\sin\theta\sin\theta_{d-1}\cos\theta_{d-2},
\\
& \cdots\cdots\cdots\cdots\cdots\cdots\cdots\cdots\cdots
\\
& x_{d-2}=R\sin\theta\sin\theta_{d-1}\cdots\cos\theta_2, \\
& x_{d-1}=R\sin\theta\sin\theta_{d-1}\cdots\sin\theta_2\cos\phi,
\\
& x_d=R\sin\theta\sin\theta_{d-1}\cdots\sin\theta_2\sin\phi,
\end{split}
\label{eq:polar}
\end{gather}
where $\phi\in[-\pi,\pi)$, $\theta_2,\ldots,\theta_d\in [0,\pi]$, and $R>0$.
On ${\mathbf S}_R^d$, the variables $\phi:=\theta_1$, $\theta:=\theta_d$, are respectively referred to as the azimuthal and radial
coordinates.

It is necessary to def\/ine the geodesic distance on $\mathbf{S}_R^d$.
In any geodesic polar coordinate system, the geodesic distance between two points $d:\mathbf{S}_R^d\times\mathbf{S}_R^d\to[0,\pi R]$ is
given in Vilenkin's polyspherical coordinates by~\cite[equation~(2)]{Cohlhypersphere}
\begin{gather}
\label{eq:cosd}
d(\mathbf{x},\mathbf{x}'):= R\cos^{-1} (\widehat{\bf x},{\widehat{\bf x}^\prime})=
R\cos^{-1}(\cos\theta\cos\theta'+\sin\theta\sin\theta'\cos\gamma),
\end{gather}
where $(\cdot,\cdot):\mathbf{R}^{d+1} \times \mathbf{R}^{d+1} \to \mathbf{R}$ is the Euclidean inner product def\/ined by
$(\mathbf{x},\mathbf{x}'):=x_0x_0'+x_1x_1'+\dots+x_dx_d'$, and~$\gamma$ is the separation angle between the two unit vectors
$\widehat{\bf x},{\widehat{\bf x}^\prime}\in{\mathbf S}^d:={\mathbf S}_1^d$.
As measured from the origin of ${\mathbf R}^{d+1}$, these unit vectors are given by $\widehat{\bf x}:={\bf x}/R$, ${\widehat{\bf
x}^\prime}:={{\bf x}^\prime}/R$.
The separation angle~$\gamma$ is def\/ined through the expression $\cos\gamma:=(\widehat{\bf x},{\widehat{\bf x}^\prime})$, and is given
in standard hyperspherical coordinates as
\begin{gather}
\label{eq:gamma}
\cos\gamma=\cos(\phi-\phi')\prod\limits_{i=1}^{d-2} \sin\theta_i\sin\theta_i' + \sum\limits_{i=1}^{d-2}
\cos\theta_i\cos\theta_i' \prod\limits_{j=1}^{i-1} \sin\theta_j\sin\theta_j'.
\end{gather}

\subsection[Fundamental solution of Laplace's equation on ${\mathbf S}_R^d$]{Fundamental solution
of Laplace's equation on $\boldsymbol{{\mathbf S}_R^d}$}

With $u\in C^2(M_d)$, where $M_d$ is a~$d$-dimensional (pseudo-)Riemannian manifold, we refer to
\begin{gather*}
\Delta u=0,
\end{gather*}
where~$\Delta$ is the Laplace--Beltrami operator on $M_d$, as Laplace's equation.
We now summarize the presentation for a~fundamental solution of Laplace's equation on ${\mathbf S}_R^d$,
given in Cohl~\cite[Theorem~3.2]{Cohlhypersphere}.
Let $\mathbf{x},\mathbf{x}' \in \mathbf{S}_R^d$, and let $\mathbf{\widehat{x}}=\mathbf{x}/R$ and $\mathbf{\widehat{x}}'=\mathbf{x}'/R$
be their unit vectors.
We f\/irst def\/ine $\mathcal{J}_d:(0,\pi)\rightarrow \mathbf{R}$ as
\begin{gather*}
\mathcal{J}_d(\Theta):= \int_{\Theta}^{\pi/2} \frac{dx}{\sin^{d-1}x},
\end{gather*}
and $\mathcal{G}_R^d:(\mathbf{S}_R^d\times\mathbf{S}_R^d) \setminus \{(\mathbf{x},\mathbf{x}):\mathbf{x} \in \mathbf{S}_R^d\}
\rightarrow \mathbf{R}$ as
\begin{gather}
\mathcal{G}_R^d(\mathbf{x},\mathbf{x}'):=\frac{\Gamma(d/2)}{2\pi^{d/2}R^{d-2}} \mathcal{J}_d(\Theta),
\label{fundsoldefnJ}
\end{gather}
where $\Theta:=\cos^{-1}((\mathbf{\widehat{x}},\mathbf{\widehat{x}}'))$ is the geodesic distance between $\mathbf{\widehat{x}}$ and
$\mathbf{\widehat{x}}'$.
Then $\mathcal{G}_R^d$ is a~fundamental solution for $-\Delta$, where~$\Delta$ is the Laplace--Beltrami operator
on $\mathbf{S}_R^d$ given~by
\begin{gather}
\label{eq:delf}
\Delta f=\frac{1}{R^2} \left[\frac{\partial^2f}{\partial\Theta^2} + (d-1) \cot\Theta \frac{\partial f}{\partial\Theta} +
\frac{1}{\sin^2\Theta} \Delta_{\mathbf{S}^{d-1}} f \right],
\end{gather}
where $\Delta_{\mathbf{S}^{d-1}}$ is the corresponding Laplace--Beltrami operator on $\mathbf{S}^{d-1}$.
Moreover,
\begin{gather}
\label{eq:Jd}
\mathcal{J}_d(\Theta)=\frac{(d-2)!}{\Gamma(d/2)2^{d/2-1}}\frac{1}{\sin^{d/2-1}\Theta} \mathsf{Q}_{d/2-1}^{1-d/2} (\cos\Theta),
\end{gather}
where $\mathsf{Q}_{\nu}^{-\nu}:(-1,1)\to{\mathbf C}$ is a~Ferrers function of the second kind def\/ined as~\cite[equation~(30)]{Cohlhypersphere}
\begin{gather*}
\mathsf{Q}_{\nu}^{-\nu}(x):=\frac{\sqrt{\pi} x(1-x^2)^{\nu/2}}{2^{\nu} \Gamma\left(\nu+\frac{1}{2}\right)}\, {}_2F_1  \left(
\begin{array}{@{}c@{}} \frac{1}{2},\nu+1 \\
\frac{3}{2} \end{array}; x^2
\right),
\end{gather*}
where $\nu\in{\mathbf C}\setminus\{-\frac12,-\frac32,-\frac52,\ldots\}$ and ${}_2F_1$ is the Gauss hypergeometric function
(see~\cite[Chap\-ter~15]{NIST}).
On Euclidean space $\mathbf{R}^d$, a~fundamental solution for $-\Delta$, the Newtonian potential, $\mathcal{N}^d:\mathbf{R}^d \times
\mathbf{R}^d \setminus \{(\mathbf{x},\mathbf{x}): \mathbf{x} \in \mathbf{R}^d \}$ is given as follows (see for
example~\cite[p.~75]{Fol4})
\begin{gather}
\label{eq:bolG}
\mathcal{N}^d(\mathbf{x},\mathbf{x}'):=
\begin{cases}
\dfrac{\Gamma(d/2)}{2\pi^{d/2} (d-2)} \Vert\mathbf{x}-\mathbf{x}'\Vert^{2-d} & \text{if}\quad d=1 \quad\text{or}\quad d \geq 3,
\\
\dfrac{1}{2\pi} \log \Vert\mathbf{x}-\mathbf{x}'\Vert^{-1} & \text{if}\quad d=2.
\end{cases}
\end{gather}

\subsection[Harmonics on the~$d$-dimensional hypersphere]{Harmonics on the~$\boldsymbol{d}$-dimensional hypersphere}

The angular harmonics are eigenfunctions of the Laplace--Beltrami operator on $\mathbf{S}^{d-1}$ which satisfy the following eigenvalue problem:
\begin{gather}
\label{eq:DELY}
\Delta_{\mathbf{S}^{d-1}} Y_l^K (\mathbf{\widehat{x}})=-l(l+d-2) Y_l^K (\mathbf{\widehat{x}}),
\end{gather}
where $\mathbf{\widehat{x}} \in \mathbf{S}^{d-1}$, $Y_l^K (\mathbf{\widehat{x}})$ are normalized hyperspherical harmonics, $l \in
\mathbf{N}_0$ is the angular momentum quantum number, and~$K$ stands for the set of $(d-2)$-quantum numbers identifying degenerate
harmonics for each~$l$~\cite{CohlKalII}.

The Riemannian (volume) measure $d\vol_g$ of the Riemannian manifold $\mathbf{S}_R^d$ with Riemannian metric~$g$ in standard
hyperspherical coordinates~\eqref{eq:polar} is given~by
\begin{gather}
\label{eq:vol}
d\vol_g=R^d\sin^{d-1}\theta d\theta d\omega:= R^d\sin^{d-1}\theta \sin^{d-2}\theta_{d-1} \cdots \sin\theta_2 d\theta
d\theta_{d-1} \cdots d\theta_2 d\phi.
\end{gather}
In computing a~fundamental solution of Laplace's equation on $\mathbf{S}_R^d$, we know that
\begin{gather}
\label{eq:negdel}
-\Delta\mathcal{G}_R^d(\mathbf{x},\mathbf{x}')=\delta_g(\mathbf{x},\mathbf{x}'),
\end{gather}
where $\delta_g(\mathbf{x},\mathbf{x}')$ is the Dirac delta distribution on the diagonal manifold $\mathbf{S}_R^d\times\mathbf{S}_R^d$.
This is the Schwartz kernel of the identity map (thought of as a~scalar distribution.) The Dirac delta distribution on the Riemannian
manifold $\mathbf{S}_R^d$ with metric~$g$ for an open set $U \subset \mathbf{S}_R^d$ with $\mathbf{x},\mathbf{x}' \in \mathbf{S}_R^d$
satisf\/ies
\begin{gather}
\label{integraldirac}
\int_U \delta_g(\mathbf{x},\mathbf{x}') d\vol_g=
\begin{cases}
1 & \text{if}\quad \mathbf{x}' \in U,
\\
0 & \text{if}\quad \mathbf{x}' \notin U.
\end{cases}
\end{gather}
Using~\eqref{eq:vol} and~\eqref{integraldirac} in standard hyperspherical coordinates~\eqref{eq:polar} on $\mathbf{S}_R^d$, the Dirac
delta distribution is given~by
\begin{gather}
\label{eq:delta}
\delta_g(\mathbf{x},\mathbf{x}')=\frac{\delta(\theta-\theta')\delta(\phi-\phi')
\delta(\theta_2-\theta_2')\cdots\delta(\theta_{d-2}-\theta_{d-2}')} {R^d\sin^{d-1}\theta'\sin\theta_2'\cdots\sin^{d-2}\theta_{d-1}'}.
\end{gather}

Due to the fact that the space ${\mathbf S}_R^d$ is homogeneous with respect to its isometry group \mbox{${\rm SO}(d+1)$}, it is therefore an
isotropic manifold.
Therefore there exists a~fundamental solution of the Laplace--Beltrami operator with purely radial, spherically symmetric dependence.
This will be important when we perform our analysis in Section~\ref{sec:gegen}.
The radial part of Laplace's equation on~$\mathbf{S}_R^d$ satisf\/ies the following dif\/ferential equation
\begin{gather*}
\frac{d^2u}{d\theta^2}+(d-1)\cot\theta \frac{du}{d\theta}-\frac{l(l+d-2)}{\sin^2\theta} u=0.
\end{gather*}
The above dif\/ferential equation is a~(transformed) associated Legendre dif\/ferential equation and originates from separating out the
radial part of~\eqref{eq:delf}.
Four solutions of this ordinary dif\/ferential equation $u_{1\pm}^{d,l},u_{2\pm}^{d,l}:(-1,1) \to \mathbf{C}$ are given
by~\cite{CohlKalII}
\begin{gather}
\label{eq:u1}
u_{1\pm}^{d,l} (\cos\theta):=\frac{1}{\sin^{d/2-1}\theta} \mathsf{P}_{d/2-1}^{\pm(d/2-1+l)} (\cos\theta),
\\
\label{eq:u2}
u_{2\pm}^{d,l} (\cos\theta):=\frac{1}{\sin^{d/2-1}\theta} \mathsf{Q}_{d/2-1}^{\pm(d/2-1+l)} (\cos\theta),
\end{gather}
where $\mathsf{P}_{\nu}^{\mu},\mathsf{Q}_{\nu}^{\mu}: (-1,1) \to \mathbf{C}$ are Ferrers functions of the f\/irst and second kind,
respectively~\cite[Section~14.3(i)]{NIST}.

\subsection{The f\/lat-space limit from hyperspherical to Euclidean geometry}\label{Theflatspacelimit}

Because each point on a~manifold is locally Euclidean, one expects that quantities def\/ined in hyperspherical geometry should reduce to
their f\/lat-space counterparts in the f\/lat-space limit.
The f\/lat-space limit is a~procedure which allows one to compare mathematical quantities def\/ined on on~$d$-dimensional spaces of constant
curvature (such as hyperspherical and hyperbolic geometry) with corresponding quantities def\/ined on~$d$-dimensional Euclidean space.
In hyperspherical geo\-metry, this method has previously been used in many dif\/ferent
contexts~\cite{Barutetal87,GroscheKarPogSis97,GroschePogSis95,HakobyanPogSisVin99}.
The limi\-ting case in 3-dimensions from the hypersphere to Euclidean space, are treated in all Helmholtz separable coordinate systems
in~\cite{PogYakh2009}.
These limiting processes are often termed (In\"on\"u--Wigner) contractions.
For contractions of Helmholtz subgroup-type separable coordinate systems from ${\mathbf S}_R^d \to {\mathbf R}^d$, see~\cite{IPSWa}.

To obtain the corresponding Euclidean limit for a~quantity expressed in a~specif\/ic coordinate system on the curved space, one performs
an asymptotic expansion as $R\to\infty$ of that quantity with coordinates approaching zero.
Take for instance the example given just below~\cite[equation~(5)]{Cohlhypersphere}.
In standard hyperspherical coordinates~\eqref{eq:polar} on ${\mathbf S}_R^d$,
the Riemannian volume measure is given by~\eqref{eq:vol}.
The f\/lat-space limit of the Riemannian volume measure on ${\mathbf S}_R^d$ is applied by taking $R\to\infty$ with the geodesic radial
distance $r=R\Theta$ is f\/ixed as measured from the origin $\theta=0$.
In the f\/lat-space limit
\begin{gather*}
R^d\sin^{d-1}\theta d\theta\sim R^d\sin^{d-1}\left(\frac{r}{R}\right)R^{-1}dr \sim R^d\left(\frac{r}{R}\right)^{d-1}R^{-1}dr=r^{d-1}dr,
\end{gather*}
which is exactly the radial volume element in standard hyperspherical coordinates on~${\mathbf R}^d$.
This example demonstrates how the f\/lat-space limit can be used (and has been used extensively) to show how Laplace--Beltrami separable
coordinate systems on ${\mathbf S}_R^d$ (such as~\eqref{eq:polar}) reduce to Laplace separable coordinate systems (such as polar or
standard hyperspherical coordinates) on~${\mathbf R}^d$.

Another interesting example of the f\/lat-space limit which we will use in Section~\ref{ApplicationsinPotentialTheory}, is Hopf
coordinates (sometimes referred to as cylindrical coordinates) on ${\mathbf S}_R^3$.
Using Hopf coordinates, points on ${\mathbf S}_R^3$ are parametrized using
\begin{gather}
x_0=R\cos\vartheta\cos\phi_1,
\qquad
x_1=R\cos\vartheta\sin\phi_1,
\nonumber\\
x_2=R\sin\vartheta\cos\phi_2,
\qquad
x_3=R\sin\vartheta\sin\phi_2,
\label{eq:hopf}
\end{gather}
where $\vartheta\in[0,\frac{\pi}{2}]$, and $\phi_1,\phi_2\in[-\pi,\pi)$.
Note that $({\bf x},{\bf x})=R^2$.
In the f\/lat-space limit we f\/ix $r=R\vartheta$, $z=R\phi_1$ and let $R\to\infty$.
In the f\/lat-space limit
\begin{gather*}
{\bf x} \sim \left(R\cos\left(\frac{r}{R}\right)\cos\left(\frac{z}{R}\right), R\cos\left(\frac{r}{R}\right)\sin\left(\frac{z}{R}\right),
R\sin\left(\frac{r}{R}\right)\cos\phi_2, R\sin\left(\frac{r}{R}\right)\sin\phi_2\right)
\\
\phantom{{\bf x}}
 \sim \left(R, R\left(\frac{z}{R}\right), R\left(\frac{r}{R}\right)\cos\phi_2, R\left(\frac{r}{R}\right)\sin\phi_2\right)
 \sim (R,z, r\cos\phi_2, r\sin\phi_2),
\end{gather*}
with the f\/irst coordinate approaching inf\/inity and the last three coordinates being circular cylindrical coordinates on ${\mathbf R}^3$.

\section{Fourier expansions for a~fundamental solution\\ of Laplace's equation}

A fundamental solution of Laplace's equation on the unit radius hypersphere $\mathfrak{g}^d:(\mathbf{S}^d \times \mathbf{S}^d) \setminus
\{(\mathbf{x},\mathbf{x}):\mathbf{x} \in \mathbf{S}^d\} \rightarrow \mathbf{R}$, is given, apart from a~multiplicative constant,
as~\eqref{fundsoldefnJ}
\begin{gather}
\label{eq:sd1}
\mathfrak{g}^d(\mathbf{\widehat{x}},\mathbf{\widehat{x}}'):=\mathcal{J}_d(d(\mathbf{\widehat{x}},\mathbf{\widehat{x}}'))
=\frac{2\pi^{d/2}R^{d-2}}{\Gamma(d/2)}\mathcal{G}_R^d(\mathbf{x},\mathbf{x}').
\end{gather}
Note that the~$d$ in $\mathfrak{g}^d$ (as in $\mathcal{G}^d$ and ${\mathbf S}^d$) is an index for $\mathfrak{g}$ representing the
dimension of the manifold and not some power of $\mathfrak{g}$.
In standard hyperspherical coordinates~\eqref{eq:polar}, we would like to perform an azimuthal Fourier expansion for a~fundamental
solution of Laplace's equation on the hypersphere; more precisely, we would like
\begin{gather}
\mathfrak{g}^d(\mathbf{\widehat{x}},\mathbf{\widehat{x}}')=\sum\limits_{m=0}^{\infty} \cos(m(\phi-\phi'))
\mathsf{G}_m^{d/2-1} (\theta,\theta',\theta_2,\theta_2',\ldots,\theta_{d-1},\theta_{d-1}'),
\label{FouriercoeffSRd}
\end{gather}
where $\mathsf{G}_m^{d/2-1}:[0,\pi]^{2d-2} \rightarrow \mathbf{R}$ which are given by
\begin{gather}
\label{eq:Smd}
\mathsf{G}_m^{d/2-1} (\theta,\theta',\theta_2,\theta_2',\ldots,\theta_{d-1},\theta_{d-1}'):=\frac{\epsilon_m}{\pi}
\int_0^{\pi} \mathfrak{g}^d(\mathbf{\widehat{x}},\mathbf{\widehat{x}}') \cos(m(\phi-\phi')) d(\phi-\phi'),
\end{gather}
where $\epsilon_m$ is the Neumann factor def\/ined as $\epsilon_m:=2-\delta_m^0$.
Using~\eqref{eq:Jd},~\eqref{eq:sd1}, we may write $\mathfrak{g}^d(\mathbf{\widehat{x}},\mathbf{\widehat{x}}')$ as
\begin{gather}
\label{eq:sd2}
\mathfrak{g}^d(\mathbf{\widehat{x}},\mathbf{\widehat{x}}')=\frac{(d-2)!}{\Gamma\left(\frac{d}{2} \right) 2^{d/2-1}(\sin
d(\mathbf{\widehat{x}},\mathbf{\widehat{x}}'))^{d/2-1}} \mathsf{Q}_{d/2-1}^{1-d/2} (\cos d(\mathbf{\widehat{x}},\mathbf{\widehat{x}}')).
\end{gather}
Through~\eqref{eq:cosd},~\eqref{eq:Smd},~\eqref{eq:sd2} in standard hyperspherical coordinates~\eqref{eq:polar}, the azimuthal Fourier
coef\/f\/icients are given~by
\begin{gather*}
\mathsf{G}_m^{d/2-1} (\theta,\theta',\theta_2,\theta_2',\ldots,\theta_{d-1},\theta_{d-1}')={
\frac{\epsilon_m(d-2)!}{\pi \Gamma(\frac{d}{2} ) 2^{d/2-1}} \int_{0}^{\pi} \frac{\mathsf{Q}_{d/2-1}^{1-d/2}
(B_d\cos\psi+A_d) \cos(m\psi)}{\big(\sqrt{1-(B_d\cos\psi+A_d)^2}\big)^{d/2-1}} d\psi},
\end{gather*}
where $\psi:=\phi-\phi'$ and $A_d,B_d:[0,\pi]^{2d-2} \rightarrow \mathbf{R}$ are def\/ined through~\eqref{eq:cosd},~\eqref{eq:gamma} as
\begin{gather*}
A_d(\theta,\theta',\theta_2,\theta_2',\ldots,\theta_{d-1},\theta_{d-1}'):=\cos\theta\cos\theta'+\sin\theta\sin\theta'
\sum\limits_{i=2}^{d-1} \cos\theta_i\cos\theta_i' \prod\limits_{j=1}^{i-1} \sin\theta_j\sin\theta_j',
\\
B_d(\theta,\theta',\theta_2,\theta_2',\ldots,\theta_{d-1},\theta_{d-1}'):=\sin\theta\sin\theta'
\prod\limits_{i=2}^{d-1} \sin\theta_i\sin\theta_i'.
\end{gather*}

\subsection{Fourier expansion for a~fundamental solution of Laplace's equation\\ on the 2-sphere}

On the unit sphere $\mathbf{S}^2$, a~fundamental solution of Laplace's equation is given by $\mathfrak{g}^2:=2\pi\mathcal{G}^2$
(cf.~\eqref{eq:sd1}), where~\cite[Theorem~3.2]{Cohlhypersphere}
\begin{gather}
\mathfrak{g}^2(\mathbf{\widehat{x}},\mathbf{\widehat{x}}') =\log\cot\frac{d(\mathbf{\widehat{x}},\mathbf{\widehat{x}}')}{2}
=\frac{1}{2}\log\frac{1+\cos d(\mathbf{\widehat{x}},\mathbf{\widehat{x}}')}{1-\cos d(\mathbf{\widehat{x}},\mathbf{\widehat{x}}')}.
\label{fs2}
\end{gather}
Note that, as in~\eqref{eq:sd1}, the 2 in $\mathfrak{g}^2$ is an index representing the dimension of the manifold and not some power of
$\mathfrak{g}$.
We can further deconstruct $\mathfrak{g}^2$ in standard hyperspherical coordinates~\eqref{eq:polar}.
Using~\eqref{eq:cosd},~\eqref{eq:gamma} to obtain
\begin{gather*}
\cos d(\mathbf{\widehat{x}},\mathbf{\widehat{x}}')=\cos\theta\cos\theta'+\sin\theta\sin\theta'\cos(\phi-\phi'),
\end{gather*}
from~\eqref{fs2} one obtains
\begin{gather}
\label{fundsol2expression}
\mathfrak{g}^2(\mathbf{\widehat{x}},\mathbf{\widehat{x}}')
=\frac{1}{2}\log\frac{1+\cos\theta\cos\theta'+\sin\theta\sin\theta'\cos(\phi-\phi')}
{1-\cos\theta\cos\theta'-\sin\theta\sin\theta'\cos(\phi-\phi')}.
\end{gather}

\begin{theorem}
Let $\phi,\phi'\in[-\pi,\pi)$. Then, the azimuthal Fourier expansion for a~fundamental solution of Laplace's equation on the unit sphere ${\mathbf S}^2$ expressed in
standard spherical coordinates~\eqref{eq:polar} is given~by
\begin{gather*}
\mathfrak{g}^2(\mathbf{\widehat{x}},\mathbf{\widehat{x}}')=\frac{1}{2} \log\frac{1+\cos\theta_{>}}{1-\cos\theta_{>}}
+ \sum\limits_{n=1}^{\infty} \frac{\cos[n(\phi-\phi')]}{n} \left(\frac{1-\cos\theta_{<}}{1+\cos\theta_{<}} \right)^{n/2}
\\
\phantom{\mathfrak{g}^2(\mathbf{\widehat{x}},\mathbf{\widehat{x}}')=}
\times \left[\left(\frac{1+\cos\theta_{>}}{1-\cos\theta_{>}} \right)^{n/2} -(-1)^n
\left(\frac{1-\cos\theta_{>}}{1+\cos\theta_{>}} \right)^{n/2} \right].
\end{gather*}
Note that we are using the notation~\eqref{alessgtr}. \label{azimuthalFourierexpansionS2}
\end{theorem}

\begin{proof}
Replacing $\psi:=\phi-\phi'$ in~\eqref{fundsol2expression} and rearranging the logarithms yields
\begin{gather}
\mathfrak{g}^2(\mathbf{\widehat{x}},\mathbf{\widehat{x}}')=\frac{1}{2} \left[\log\frac{1+\cos\theta\cos\theta'}
{1-\cos\theta\cos\theta'}+\log(1+z_{+}\cos\psi)-\log(1-z_{-}\cos\psi) \right],
\label{fundsollogexp2}
\end{gather}
where
\begin{gather}
\label{eq:z}
z_{\pm}:=\frac{\sin\theta\sin\theta'}{1\pm\cos\theta\cos\theta'}.
\end{gather}
Note that that $z_{\pm} \in (-1,1)$ for $\theta,\theta' \in [0,\pi]$.
For $x \in (-1,1)$, we have the Maclaurin series
\begin{gather*}
\log(1\pm x)=-\sum\limits_{n=1}^{\infty} \frac{(\mp 1)^nx^n}{n}.
\end{gather*}
Therefore, away from the singularity at $\mathbf{x}=\mathbf{x}'$, we have $\lambda_{\pm}: (-1,1)\times(-\pi,\pi) \to (-\infty,\log2)$
def\/ined~by
\begin{gather}
\label{eq:lpm1}
\lambda_{\pm} (z_{\pm},\psi):= \log(1\pm z_{\pm}\cos\psi)=-\sum\limits_{k=1}^{\infty}
\frac{(\mp1)^kz_{\pm}^k\cos^k(\psi)}{k}.
\end{gather}
We can expand powers of cosine using the following trigonometric identity:
\begin{gather*}
\cos^k\psi=\frac{1}{2^k} \sum\limits_{n=0}^{k} {k \choose n}\cos[(2n-k)\psi],
\end{gather*}
which is the standard expansion for powers~\cite[p.~52]{FoxParker} when using Chebyshev polynomials of the f\/irst kind $T_m(\cos\psi)=\cos(m\psi)$.
Inserting the above expression into~\eqref{eq:lpm1}, we obtain the double-summation expression
\begin{gather}
\label{eq:lpm2}
\lambda_{\pm}(z_{\pm},\psi)=-\sum\limits_{k=1}^{\infty} \sum\limits_{n=0}^{k} \frac{(\mp1)^kz_{\pm}^k}{2^kk}
{k \choose n} \cos[(2n-k)\psi].
\end{gather}

Our current goal is to identify closed-form expressions for the azimuthal Fourier cosine coef\/f\/icients $c_{n'}$ of the sum given~by
\begin{gather*}
\lambda_\pm(z_\pm,\psi)=\sum\limits_{n'=0}^\infty c_{n'}(z_\pm)\cos(n'\psi).
\end{gather*}
We accomplish this by performing a~double-index replacement as in~\cite[equation~(41)]{CohlKalII} to~\eqref{eq:lpm2}.
The double sum over~$k$,~$n$ is rearranged into two separate double sums with dif\/ferent indices~$k'$,~$n'$, f\/irst for~$k\le 2n$ and second for
$k\ge 2n$.
If $k\le 2n$, let $k'=k-n$, $n'=2n-k$ which implies $k=2k'+n'$ and $n=n'+k'$.
If $k\ge 2n$, let $k'=n$ and $n'=k-2n$, which implies $k=2k'+n'$ and $n=k'$.
Note that the identity $\binom{n}{k}=\binom{n}{n-k}$ implies $\binom{2n'+k'}{n'+k'}=\binom{2k'+n'}{k'}$.
After halving the the double counting when $k=2n$ and replacing $k'\mapsto k$ and $n'\mapsto n$, we obtain
\begin{gather*}
\lambda_{\pm}(z_{\pm},\psi)=-\sum\limits_{k=1}^{\infty} \frac{z_{\pm}^{2k}}{2^{2k}(2k)} {2k \choose k}
-2\sum\limits_{n=1}^{\infty} (\mp1)^n\cos(n\psi) \sum\limits_{k=0}^{\infty}
\frac{z_{\pm}^{2k+n}}{2^{2k+n}(2k+n)} {2k+n \choose k}.
\end{gather*}
If we def\/ine $I_{\pm}:(-1,1) \to (-\log2,0)$ and $J_{n,\pm}:(-1,1) \to (0,(\mp2)^{n}/n)$ such that
\begin{gather*}
\lambda_{\pm}(z_{\pm},\psi)=I_{\pm}(z_{\pm}) + \sum\limits_{n=0}^{\infty} J_{n,\pm}(z_{\pm}) \cos(n\psi),
\end{gather*}
then these terms are given by~\cite[Section~4.1]{CohlKalII}
\begin{gather*}
I_{\pm}(z_{\pm}):= -\sum\limits_{k=1}^{\infty} \frac{z_{\pm}^{2k}}{2^{2k}(2k)} {2k \choose k}=-\log2+\log\left(1+\sqrt{1-z_{\pm}}\right),
\\
J_{n,\pm}(z_{\pm}):= -2(\mp1)^n\sum\limits_{k=0}^{\infty} \frac{z_{\pm}^{2k+n}}{2^{2k+n}(2k+n)} {2k+n \choose k}=
\frac{2(\mp1)^n}{n} \left(\frac{1-\sqrt{1-z_{\pm}^2}}{1+\sqrt{1-z_{\pm}^2}} \right)^{n/2}.
\end{gather*}
Finally, from~\eqref{eq:z}, we have
\begin{gather*}
\lambda_{\pm}= -\log2+\log\frac{(1+\cos\theta_{<})(1\pm\cos\theta_{>})} {1\pm\cos\theta\cos\theta'}
\\
\phantom{\lambda_{\pm}=}
-\sum\limits_{n=1}^{\infty} \frac{2(\mp1)^n}{n} \cos(n\psi)
\left(\frac{(1-\cos\theta_{<})(1\mp\cos\theta_{>})}{(1+\cos\theta_{<}) (1\pm\cos\theta_{>})}\right)^{n/2}.
\end{gather*}
Utilizing this expression in~\eqref{fundsollogexp2} completes the proof.
\end{proof}

Note, as shown in~\cite[equation~(39)]{CohlKalII}, that (using the notation~\eqref{alessgtr})
\begin{gather}
\mathfrak{n}^2(\mathbf{x},\mathbf{x}')=\log \Vert\mathbf{x}-\mathbf{x}'\Vert^{-1}=-\log r_{>} +
\sum\limits_{n=1}^{\infty} \frac{\cos[n(\phi-\phi')]}{n} \left(\frac{r_{<}}{r_{>}} \right)^n,
\label{2discpotential}
\end{gather}
where $\mathfrak{n}^2:=2\pi \mathcal{N}^2$ is def\/ined using~\eqref{eq:bolG}, normalized as in~\eqref{eq:sd1}.
It is straightforward to check that $\mathfrak{g}^2(\mathbf{\widehat{x}},\mathbf{\widehat{x}}') \rightarrow
\mathfrak{n}^2(\mathbf{x},\mathbf{x}')$ as $\theta,\theta' \rightarrow 0^{+}$ in the f\/lat-space limit $R\to\infty$.

\subsection{Fourier expansion for a~fundamental solution of Laplace's equation\\ on the 3-sphere}\label{four3}

On the hypersphere $\mathbf{S}_R^3$, a~fundamental solution of Laplace's equation $\mathfrak{g}^3:=4\pi R \mathcal{G}_R^3$
(cf.~\eqref{eq:sd1}) is given by~\cite[Theorem~3.2]{Cohlhypersphere}
\begin{gather}
\mathfrak{g}^3(\mathbf{\widehat{x}},\mathbf{\widehat{x}}')=\cot(d(\mathbf{\widehat{x}},\mathbf{\widehat{x}}'))
=\frac{\cos(d(\mathbf{\widehat{x}},\mathbf{\widehat{x}}'))}{\sqrt{1-\cos(d(\mathbf{\widehat{x}},\mathbf{\widehat{x}}'))^2}},
\label{Fourierexpansiond3}
\end{gather}
where~\eqref{eq:cosd},~\eqref{eq:gamma} still hold.
As indicated in~\eqref{eq:sd1}, the 3 in $\mathfrak{g}^3$ above, should be interpreted as an index for $\mathfrak{g}$ representing the
dimension of the hypersphere, and not some power of $\mathfrak{g}$.
If we def\/ine $A:[0,\pi]^4 \to [-1,1]$ and $B:[0,\pi]^4 \to [0,1]$ respectively as{\samepage
\begin{gather}
\label{Adefn}
A(\theta,\theta',\theta_2,\theta_2') := \cos\theta\cos\theta' +\sin\theta\sin\theta'\cos\theta_2\cos\theta_2',
\\
\label{Bdefn}
B(\theta,\theta',\theta_2,\theta_2') := \sin\theta\sin\theta'\sin\theta_2\sin\theta_2',
\end{gather}

}

\noindent
then the azimuthal Fourier coef\/f\/icients $\mathsf{G}_m^{1/2}:[0,\pi]^4 \rightarrow \mathbf{R}$ of $\mathfrak{g}^3$ are expressed through
\begin{gather}
\label{g3azimuthalfouriersum}
\mathfrak{g}^3(\mathbf{\widehat{x}},\mathbf{\widehat{x}}')=\sum\limits_{m=0}^{\infty} \mathsf{G}_m^{1/2}
(\theta,\theta',\theta_2,\theta_2') \cos (m\psi),
\end{gather}
where $\psi:=\phi-\phi'$.
By application of orthogonality for Chebyshev polynomials of the f\/irst kind~\cite[Table~18.3.1]{NIST}
to~\eqref{g3azimuthalfouriersum}, using~\eqref{eq:cosd} and~\eqref{Fourierexpansiond3}, we produce an integral representation
for the azimuthal Fourier coef\/f\/icients
\begin{gather*}
\mathsf{G}_m^{1/2}(\theta,\theta',\theta_2,\theta_2')= \frac{\epsilon_m}{\pi} \int_{0}^{\pi} \frac{\left(\cos\psi+A/B
\right) \cos(m\psi) d\psi}{\sqrt{\left(\frac{1-A}{B}-\cos\psi\right) \left(\frac{1+A}{B}+\cos\psi\right)}}.
\end{gather*}
If we substitute $x=\cos\psi$, then this integral can be converted into
\begin{gather}
\label{eq:Smhalf1}
\mathsf{G}_m^{1/2} (\theta,\theta',\theta_2,\theta_2')=\frac{\epsilon_m}{\pi} \int_{-1}^{1}
\frac{(x+A/B)T_m(x)dx}{\sqrt{(1-x^2) \left(\frac{1-A}{B}-x\right)\left(\frac{1+A}{B}+x\right)}}.
\end{gather}
Since $(x+A/B)T_m(x)$ is a~polynomial in~$x$, it is suf\/f\/icient to solve the integral
\begin{gather}
\label{eq:xint}
\int_{-1}^{1} \frac{x^p dx}{\sqrt{(1-x^2)\left(\frac{1-A}{B}-x\right) \left(\frac{1+A}{B}+x\right)}},
\end{gather}
which by def\/inition is an elliptic integral since it involves the square root of a~quartic in~$x$ multiplied by a~rational function
of~$x$.

Let $F:[0,\frac{\pi}{2}] \times [0,1) \rightarrow \mathbf{R}$ be Legendre's incomplete elliptic integral of the f\/irst kind which can be
def\/ined through the def\/inite integral~\cite[Section~19.2(ii)]{NIST}
\begin{gather*}
F(\varphi,k):=\int_{0}^{\varphi} \frac{d\theta}{\sqrt{1-k^2\sin^2\theta}},
\end{gather*}
$K:[0,1) \rightarrow [1,\infty)$ is Legendre's complete elliptic integral of the f\/irst kind given by $K(k):=F
(\frac{\pi}{2},k )$, $E:[0,1] \rightarrow [1,\frac{\pi}{2}]$ and $\Pi:[0,\infty)\setminus\{1\} \times [0,1) \rightarrow
\mathbf{R}$ are Legendre's complete elliptic integrals of the second and third kind, respectively, def\/ined as
\begin{gather*}
E(k):=\int_{0}^{\pi/2} \sqrt{1-k^2\sin^2\theta} d\theta,
\qquad
\Pi(\alpha^2,k):=\int_{0}^{\pi/2} \frac{d\theta}{\sqrt{1-k^2\sin^2\theta} (1-\alpha^2\sin^2\theta)}.
\end{gather*}
For $n\in{\mathbf N}_0$, let $(\cdot)_n:{\mathbf C}\to{\mathbf C}$, denote the Pochhammer symbol (rising factorial), which is
def\/ined by $(a)_n:=(a)(a+1)\cdots(a+n-1)$.
\begin{theorem}
\label{azimuthalFourierexp3D}
Let $A$, $B$ be defined by~\eqref{Adefn},~\eqref{Bdefn} and $\alpha$, $k$ be defined by~\eqref{alphagphietcdefn}.
The azimuthal Fourier coefficients for a~fundamental solution of Laplace's equation on the unit hypersphere ${\mathbf S}^3$ expressed in
standard hyperspherical coordinates~\eqref{eq:polar} is given~by
\begin{gather*}
\mathsf{G}_m^{1/2} (\theta,\theta',\theta_2,\theta_2')=\frac{2B\epsilon_m}{\pi\sqrt{(1-A+B)(1+A-B)}}
\sum\limits_{p=0}^{m+1} \sum\limits_{j=0}^{p} a_p \frac{(B-1+A)^j}{B^p(1-A)^{j-p}} {p \choose j}
V_j(\alpha,k),
\end{gather*}
where
\begin{gather}
a_p =
\begin{cases}
\displaystyle
\frac{A(m)_m (\frac{-m}{2})_{\frac{m-p}{2}} (\frac{-m+1}{2})_{\frac{m-p}{2}}}{2^mB(\frac12)_m(1-m)_{\frac{m-p}{2}}(\frac{m-p}{2})!}
 & \text{if}\quad p=m-2\left\lfloor\frac{m}{2}\right\rfloor,\ldots,m,
\vspace{1mm}\\
\displaystyle
\frac{(m)_m (\frac{-m}{2})_{\frac{m-p+1}{2}} (\frac{-m+1}{2})_{\frac{m-p+1}{2}}}{2^m(\frac12)_m(1-m)_{\frac{m-p+1}{2}}(\frac{m-p+1}{2})!}
 & \text{if}\quad
p=m-2\left\lfloor\frac{m}{2}\right\rfloor+1,\ldots,m+1,
\end{cases}
\label{coefap}
\\
V_0(\alpha,k)=K(k), V_1(\alpha,k)=\Pi(\alpha^2,k),
\nonumber
\\
V_2(\alpha,k)=\frac{1}{2(\alpha^2-1)(k^2-\alpha^2)}\big[(k^2-\alpha^2)K(k)+\alpha^2E(k)
\nonumber
\\
\phantom{V_2(\alpha,k)=}{}
 +(2\alpha^2k^2+2\alpha^2-\alpha^4-3k^2)\Pi(\alpha^2,k)\big],
\nonumber
\end{gather}
and the values of $V_j(\alpha,k)$ for $j>3$ can be computed using the recurrence relation~{\rm \cite[equation~(336.00-03)]{ByrdFriedman}}
\begin{gather*}
V_{j+3}(\alpha,k)=\frac{1}{2(j+2)(\alpha^2-1)(k^2-\alpha^2)} \Bigl[(2j+1)k^2V_j(\alpha,k)
\\
\phantom{V_{j+3}(\alpha,k) =}{}
+2(j+1)(\alpha^2k^2+\alpha^2-3k^2)V_{j+1}(\alpha,k) \\
\phantom{V_{j+3}(\alpha,k) =}{}
+(2j+3)(\alpha^4-2\alpha^2k^2-2\alpha^2+3k^2)V_{j+1}(\alpha,k)\Bigr].
\end{gather*}
\end{theorem}

\begin{proof}
We can directly compute~\eqref{eq:Smhalf1} using~\cite[equation~(255.17)]{ByrdFriedman}. If we def\/ine
\begin{gather}
\label{eq:ineq}
a:=\frac{1-A}{B}, \qquad b:=1, \qquad c:=-1, \qquad y:=-1, \qquad d:=-\frac{1+A}{B},
\end{gather}
(clearly $d<y\le c<b<a$), then we can express the Fourier coef\/f\/icient~\eqref{eq:Smhalf1} as a~linear combination of integrals
\begin{gather}
\label{eq:elpint}
\int_{-1}^{1} \frac{x^p dx}{\sqrt{(a-x)(b-x)(x-c)(x-d)}}=\frac{2}{\sqrt{(a-c)(b-d)}} \int_{0}^{u_1}
\left(\frac{1-\alpha_1^2 {\sn}^2u}{1-\alpha^2 {\sn}^2u} \right)^p du,
\end{gather}
where $0\leq p \leq m+1$.
In this expression, ${\sn} u$ is a~Jacobi elliptic function. Byrd and Friedman~\cite{ByrdFriedman}
give a~procedure for computing~\eqref{eq:elpint} for all $p \in \mathbf{N}_0$.
These integrals will be given in terms of complete elliptic integrals of the f\/irst three kind.
We have the following def\/initions from~\cite[Section~255, equations~(255.17), (340.00)]{ByrdFriedman}:
\begin{gather*}
0< \alpha^2=\frac{b-c}{a-c}<k^2,
\qquad g=\frac{2}{\sqrt{(a-c)(b-d)}},
\qquad \varphi=\sin^{-1} \sqrt{\frac{(a-c)(b-y)}{(b-c)(a-y)}},
\\
\alpha_1^2=\frac{a(b-c)}{b(a-c)},
\qquad
k^2=\frac{(b-c)(a-d)}{(a-c)(b-d)},
\qquad u_1=F(\varphi,k).
\end{gather*}
For our specif\/ic choices in~\eqref{eq:ineq}, these reduce to
\begin{gather}
\alpha^2=\frac{2B}{1-A+B},
\qquad
g=\frac{2B}{\sqrt{(1+A+B)(1-A+B)}},
\qquad
\varphi=\frac{\pi}{2},
\nonumber
\\
\alpha_1^2=\frac{2(1-A)}{1-A+B},
\qquad
k^2=\frac{4B}{(1+A+B)(1-A+B)},
\qquad
u_1=K(k).
\label{alphagphietcdefn}
\end{gather}
The required integration formula (cf.~\cite[equation~(340.04)]{ByrdFriedman}) is given~by
\begin{gather*}
\int_{-1}^{1} \frac{x^p dx}{\sqrt{(a-x)(b-x)(x-c)(x-d)}}=\frac{g(b\alpha_1^{2})^p}{\alpha^{2p}}
\sum\limits_{j=0}^{p} \frac{(\alpha^2-\alpha_1^2)^j}{\alpha_1^{2j}} {p \choose j} V_j(\alpha,k).
\end{gather*}
With our specif\/ic values, we have{\samepage
\begin{gather}
\int_{-1}^{1} \frac{x^p dx}{\sqrt{(\frac{1-A}{B}-x)(1-x)(x+1) (x+\frac{1+A}{B})}}
\nonumber
\\
\qquad=\frac{2B}{\sqrt{(1+A+B)(1-A+B)}}\left(\frac{1-A}{B}\right)^p \sum\limits_{j=0}^{p}\left(\frac{A+B-1}{1-A} \right)^j {p \choose j} V_j(\alpha,k).
\label{sumforintegralforfouriercoef}
\end{gather}
}

Since we have the integral~\eqref{eq:xint}, we would like to express $(x+A/B)T_m(x)$ in a~sum over powers of~$x$.
This reduces to determining the coef\/f\/icients $a_p$ such that
\begin{gather}
(x+A/B) T_m(x)=\sum\limits_{p=0}^{m+1} a_px^p.
\label{xpABTmx}
\end{gather}
The Chebyshev polynomial $T_m(x)$ is expressible as a~f\/inite sum in terms of the binomial $(1-x)^m$ (see~\cite[equation~(15.9.5)]{NIST})
as $T_m(x)={}_2F_1(-m,m;\frac12;\frac{1-x}{2})$. We can further expand these powers for all $m\in{\mathbf N}_0$ to obtain
\begin{gather}
T_m(x)=\frac{(m)_m}{2^m\left(\frac12\right)_m} \sum\limits_{k=0}^{\left\lfloor\frac{m}{2}\right\rfloor}
\frac{\left(\frac{-m}{2}\right)_k \left(\frac{-m+1}{2}\right)_k}{(1-m)_kk!}x^{m-2k}.
\label{newChebyshev}
\end{gather}
(Note that $(m)_m$ can only be written in terms of factorials for $m\ge 1$.) Examining~\eqref{newChebyshev} multiplied by $x+A/B$
produces the coef\/f\/icients $a_p$ in~\eqref{xpABTmx} as~\eqref{coefap}.
Combining~\eqref{eq:Smhalf1}, ~\eqref{sumforintegralforfouriercoef}, and~\eqref{xpABTmx} completes the proof.
\end{proof}

Utilizing the above procedure to obtain the azimuthal Fourier coef\/f\/icients of a~fundamental solution of Laplace's equation on ${\mathbf
S}_R^3$ in terms of complete elliptic integrals, let us directly compute the $m=0$ component.
In this case~\eqref{eq:Smhalf1} reduces to
\begin{gather*}
\mathsf{G}_0^{1/2} (\theta,\theta',\theta_2,\theta_2')
=\frac{1}{\pi} \int_{-1}^{1}\frac{(x+A/B)dx}{\sqrt{(1-x)(1+x)(\frac{1-A}{B}-x) (\frac{1+A}{B}+x)}}.
\end{gather*}
Therefore, using the above formulas, we have
\begin{gather*}
\mathcal{G}_R^3({\bf x},{{\bf x}^\prime})\big|_{m=0}=\frac{1}{4\pi R}\mathsf{G}_0^{1/2}(\theta,\theta',\theta_2,\theta_2')
= \frac{\left(V_0(\alpha,k) + (B-1+A) V_1(\alpha,k) \right)} {2\pi^2 R\sqrt{(1-A+B)(1+A-B)}}
\\
\phantom{\mathcal{G}_R^3({\bf x},{{\bf x}^\prime})\big|_{m=0}}{}
= \frac{1}{2\pi^2 R} \left\{K(k) + [\cos\theta\cos\theta'+\sin\theta\sin\theta' \cos(\theta_2-\theta_2')-1] \Pi\big(\alpha^2,k\big)\right\}
\\
\phantom{\mathcal{G}_R^3({\bf x},{{\bf x}^\prime})\big|_{m=0}=}{}
\times [1+\cos\theta\cos\theta' +\sin\theta\sin\theta'\cos(\theta_2+\theta_2')]^{1/2}
\\
\phantom{\mathcal{G}_R^3({\bf x},{{\bf x}^\prime})\big|_{m=0}=}{}
\times [1-\cos\theta\cos\theta' -\sin\theta\sin\theta'\cos(\theta_2+\theta_2')]^{1/2}.
\end{gather*}

The Fourier expansion for a~fundamental solution of Laplace's equation in three-dimensional Euclidean space in standard spherical
coordinates $\mathbf{x}=(r\sin\theta\cos\phi, r\sin\theta\sin\phi, r\cos\theta)$ is given by~\eqref{eq:bolG} and~\cite[equation~(1.3)]{CRTB} as
\begin{gather*}
\mathfrak{n}^3(\mathbf{x},\mathbf{x}')={ \frac{1}{\Vert\mathbf{x}-\mathbf{x}'\Vert}}
\\
\hphantom{\mathfrak{n}^3(\mathbf{x},\mathbf{x}')}{}
=
\frac{1}{\pi\sqrt{rr'\sin\theta\sin\theta'}} \sum\limits_{m=-\infty}^{\infty} e^{im(\phi-\phi')} Q_{m-1/2}
\left(\frac{r^2+{r'}^2-2rr'\cos\theta\cos\theta'} {2rr'\sin\theta\sin\theta'} \right),
\end{gather*}
where $\mathfrak{n}^3:=4\pi\mathcal{N}^3$ is normalized as in~\eqref{eq:sd1}~\cite{CT}.
By~\cite[equation~(8.13.3)]{Abra},
\begin{gather*}
{\mathcal N}^3(\mathbf{x},\mathbf{x}')\big|_{m=0}=\frac{1} {2\pi^2\sqrt{r^2+{r'}^2-2rr'\cos(\theta+\theta')}} K
\left(\sqrt{\frac{4rr'\sin\theta\sin\theta'}{r^2+{r'}^2-2rr'\cos(\theta+\theta')}} \right),
\end{gather*}
which is the f\/lat-space limit of $\mathcal{G}_R^3({\bf x},{{\bf x}^\prime})\big|_{m=0}$ as $\theta,\theta' \rightarrow 0^{+}$, $R\to\infty$.

\section[Gegenbauer polynomial expansions on the~$d$-sphere]{Gegenbauer polynomial expansions on the~$\boldsymbol{d}$-sphere}\label{sec:gegen}

The Gegenbauer polynomial $C_l^{\mu}:\mathbf{C} \to \mathbf{C}$, $l \in \mathbf{N}_0$, $\mu\in(-\frac12,\infty)\setminus\{0\}$, can be
def\/ined in terms of the Gauss hypergeometric function as
\begin{gather*}
C_l^{\mu}(x):=\frac{(2\mu)_l}{l!}\, \hyp21{-l,2\mu+l}{\mu+\frac{1}{2}}{\frac{1-x}{2}}.
\end{gather*}
Note that in the following theorem, we are using the notation~\eqref{alessgtr}.

\begin{theorem}\label{Gegexpansion}
The Gegenbauer polynomial expansion for a~fundamental solution of Laplace's equation on the~$R$-radius,~$d$-dimensional hypersphere
${\mathbf S}_R^d$ for $d\ge 3$ expressed in standard hyperspherical coordinates~\eqref{eq:polar} is given~by
\begin{gather}
\mathcal{G}_R^d(\mathbf{x},\mathbf{x}')=\frac{\Gamma(d/2)}{(d-2)R^{d-2}(\sin\theta\sin\theta')^{d/2-1}}\nonumber
\\
\hphantom{\mathcal{G}_R^d(\mathbf{x},\mathbf{x}')=}{}\times\begin{cases}
\displaystyle  \frac{(-1)^{d/2-1}}{2\pi^{d/2}} \sum\limits_{l=0}^{\infty} (-1)^l(2l+d-2)
\mathsf{P}_{d/2-1}^{-(d/2-1+l)} (\cos\theta_{<}) & \\
\qquad {}\times \mathsf{Q}_{d/2-1}^{d/2-1+l} (\cos\theta_{>}) C_l^{d/2-1} (\cos\gamma)
 & \text{if}\  d \ \text{is even},
\vspace{1mm}\\
\displaystyle \frac{(-1)^{(d-3)/2}}{4\pi^{d/2-1}} \sum\limits_{l=0}^{\infty} (-1)^l(2l+d-2) \mathsf{P}_{d/2-1}^{-(d/2-1+l)}
(\cos\theta_{<}) & \\
\qquad {}\times \mathsf{P}_{d/2-1}^{d/2-1+l} (\cos\theta_{>}) C_l^{d/2-1} (\cos\gamma)
 & \text{if}\  d \ \text{is odd}.
\end{cases}\label{eq:SDY4}
\end{gather}
\end{theorem}

Theorem~\ref{azimuthalFourierexpansionS2} can be obtained by starting with Theorem~\ref{Gegexpansion} and taking the limit as $d\to 2$.

\begin{proof}
The completeness relation for hyperspherical harmonics in standard hyperspherical coordinates~\eqref{eq:polar} is given~by
\begin{gather*}
\sum\limits_{l=0}^{\infty} \sum\limits_{K} Y_l^K (\phi,\theta_2,\ldots,\theta_{d-1}) \overline{Y_l^K
(\phi',\theta_2',\ldots,\theta_{d-1}')}=\frac{\delta(\phi-\phi') \delta(\theta_2 - \theta_2') \cdots
\delta(\theta_{d-1}-\theta_{d-1}')}{\sin^{d-2}\theta_{d-1}'\cdots\sin\theta_2'}.
\end{gather*}
Therefore, through~\eqref{eq:delta}, we can write
\begin{gather}
\label{eq:deltay}
\delta_g(\mathbf{x},\mathbf{x}')=\frac{\delta(\theta-\theta')}{R^d\sin^{d-1}\theta'} \sum\limits_{l=0}^{\infty}
\sum\limits_{K} Y_l^K (\phi,\theta_2,\ldots,\theta_{d-1}) \overline{Y_l^K (\phi',\theta_2',\ldots,\theta_{d-1}')}.
\end{gather}
Since $\mathcal{G}_R^d$ is harmonic on its domain for f\/ixed $\phi'\in[-\pi,\pi)$, $\theta,\theta',\theta_2',\ldots,\theta_{d-1}' \in
[0,\pi]$, its restriction is in $C^2(\mathbf{S}^{d-1})$, and therefore has a~unique expansion in hyperspherical harmonics, namely
\begin{gather}
\label{eq:SDY1}
\mathcal{G}_R^d(\mathbf{x},\mathbf{x}')= \sum\limits_{l=0}^{\infty} \sum\limits_{K} u_l^K
(\theta,\theta',\phi',\theta_2',\ldots,\theta_{d-1}') Y_l^K (\phi,\theta_2,\ldots,\theta_{d-1}),
\end{gather}
where $u_l^K:[0,\pi]^2\times[-\pi,\pi)\times[0,\pi]^{d-2} \to \mathbf{C}$.
If we substitute~\eqref{eq:deltay},~\eqref{eq:SDY1} into~\eqref{eq:negdel} and use~\eqref{eq:delf},~\eqref{eq:DELY}, we obtain
\begin{gather}
- \sum\limits_{l=0}^{\infty} \sum\limits_K Y_l^K (\phi,\theta_2,\ldots,\theta_{d-1})
\left[\frac{\partial^2}{\partial\theta^2} +(d-1)\cot\theta\frac{\partial}{\partial\theta} -\frac{l(l+d-2)}{\sin^2\theta} \right]\nonumber\\
\qquad\quad{}\times u_l^K
(\theta,\theta',\phi',\theta_2',\ldots,\theta_{d-1}').
\nonumber
\\
\label{eq:SUMY1}
\qquad{}
= \sum\limits_{l=0}^{\infty} \sum\limits_{K} Y_l^K (\phi,\theta_2,\ldots,\theta_{d-1}) \overline{Y_l^K
(\phi',\theta'_2,\ldots,\theta_{d-1}')} \frac{\delta(\theta-\theta')}{R^d\sin^{d-1}\theta'}.
\end{gather}
This indicates the ansatz $u_l:[0,\pi]^2 \to \mathbf{R}$ given~by
\begin{gather}
\label{eq:ul1}
u_l^K (\theta,\theta',\phi',\theta_2',\ldots,\theta_{d-1}')=u_l(\theta,\theta') \overline{Y_l^K
(\phi',\theta_2',\ldots,\theta_{d-1}')}.
\end{gather}
What is derived subsequently is an explicit form of a~fundamental solution for the associated Legendre equation.
From~\eqref{eq:SDY1},~\eqref{eq:ul1}, the expression for a~fundamental solution of the Laplace--Beltrami operator in standard
hyperspherical coordinates on the hypersphere is given~by
\begin{gather}
\label{eq:SDY2}
\mathcal{G}_R^d(\mathbf{x},\mathbf{x}')= \sum\limits_{l=0}^{\infty} u_l(\theta,\theta') \sum\limits_{K} Y_l^K
(\phi,\theta_2,\ldots,\theta_{d-1}) \overline{Y_l^K (\phi',\theta_2,\ldots,\theta_{d-1}')}.
\end{gather}
We use the addition theorem for hyperspherical harmonics~\cite{WenAvery} given~by
\begin{gather}
\label{eq:SUMY2}
\sum\limits_{K} Y_l^K (\mathbf{\widehat{x}}) \overline{Y_l^K (\mathbf{\widehat{x}}')}=
\frac{\Gamma(d/2)}{2\pi^{d/2}(d-2)}(2l+d-2) C_l^{d/2-1} (\cos\gamma),
\end{gather}
where~$\gamma$ is the separation angle (see Section~\ref{Coordinates}).
The above equation~\eqref{eq:SDY2} can now be simplif\/ied using~\eqref{eq:SUMY2}; therefore,
\begin{gather}
\label{eq:SDY3}
\mathcal{G}_R^d(\mathbf{x},\mathbf{x}')= \frac{\Gamma(d/2)}{2\pi^{d/2}(d-2)} \sum\limits_{l=0}^{\infty} u_l(\theta,\theta')
(2l+d-2) C_l^{d/2-1} (\cos\gamma).
\end{gather}

Now we compute the exact expression of $u_l(\theta,\theta')$.
By separating out the radial part in~\eqref{eq:SUMY1} and using~\eqref{eq:ul1}, we obtain the dif\/ferential equation
\begin{gather}
\label{eq:ul2}
\frac{\partial^2u_l(\theta,\theta')}{\partial\theta^2} + (d-1)\cot\theta\frac{\partial u_l(\theta,\theta')}{\partial\theta} -
\frac{l(l+d-2)u_l(\theta,\theta')}{\sin^2\theta}=-\frac{\delta(\theta-\theta')}{R^d\sin^{d-1}\theta'}.
\end{gather}
Away from $\theta=\theta'$, solutions to the dif\/ferential equation~\eqref{eq:ul2} must be given by linearly independent
solutions to the homogeneous equation, which are given by~\eqref{eq:u1},~\eqref{eq:u2}.
These solutions are given in terms of Ferrers functions of the f\/irst and second kind.
The selection of linear independent pairs of Ferrers functions depends on the dimension of the space.
We make use of the Wronskian relations~\cite[equation~(14.2.6)]{NIST},~\cite[p.~170]{MOS}, to obtain
\begin{gather}
\label{WPQ}
W\left\{\mathsf{P}_{d/2-1}^{-(d/2-1+l)}(\cos\theta'), \mathsf{Q}_{d/2-1}^{d/2-1+l}(\cos\theta') \right\} =
\frac{(-1)^{d/2-1+l}}{\sin^2\theta'},
\\[0.2cm]
\label{WPP}
W\left\{\mathsf{P}_{d/2-1}^{-(d/2-1+l)}(\cos\theta'), \mathsf{P}_{d/2-1}^{d/2-1+l}(\cos\theta') \right\} =
\frac{2(-1)^{(d-3)/2+l}}{\pi\sin^2\theta'},
\end{gather}
for~$d$ even and~$d$ odd respectively.
The Ferrers functions provide a~linearly independent set of solutions for the homogeneous problem provided that the Wronskians are
non-zero and well-def\/ined.
The f\/irst Wronskian relation~\eqref{WPQ} is well-def\/ined for~$d$ even, but is not def\/ined for~$d$ odd.
The second Wronskian relation is well-def\/ined for~$d$ odd, but is not def\/ined for~$d$ even.
So unlike the corresponding problem in hyperbolic geometry~\cite[equation~(70)]{CohlKalII}, on the hypersphere we must choose our
linearly dependent solutions dif\/ferently depending on the evenness or oddness of the dimension.
Therefore, the solution to~\eqref{eq:ul2} is given~by
\begin{gather}
\label{eq:ul3}
u_l(\theta,\theta')=
\begin{cases}
\displaystyle
\frac{C}{(\sin\theta\sin\theta')^{d/2-1}}\mathsf{P}_{d/2-1}^{-(d/2-1+l)} (\cos\theta_{<})\mathsf{Q}_{d/2-1}^{d/2-1+l}(\cos\theta_{>})
 & \text{if}\  d \ \text{is even},
\vspace{1mm}\\
\displaystyle
\frac{D}{(\sin\theta\sin\theta')^{d/2-1}}\mathsf{P}_{d/2-1}^{-(d/2-1+l)} (\cos\theta_{<}) \mathsf{P}_{d/2-1}^{d/2-1+l}(\cos\theta_{>})
 & \text{if}\  d \ \text{is odd},
\end{cases}
\end{gather}
such that $u_l(\theta,\theta')$ is continuous at $\theta=\theta'$ and $C,D \in \mathbf{R}$.

In order to determine the constants $C$, $D$, we f\/irst def\/ine
\begin{gather}
\label{eq:vl}
v_l(\theta,\theta'):= (\sin\theta\sin\theta')^{(d-1)/2} u_l(\theta,\theta').
\end{gather}
Thus~\eqref{eq:ul2} can be given~by
\begin{gather*}
\frac{\partial^2v_l(\theta,\theta')}{\partial\theta^2} - \frac{v_l(\theta,\theta')}{4}\left[-5(d-1)^2 +
\frac{5d^2-12d+7+4l^2+4ld-8l}{\sin^2\theta} \right]=-\frac{\delta(\theta-\theta')}{R^{d-2}},
\end{gather*}
which we then integrate over~$\theta$ from $\theta'-\epsilon$ to $\theta'+\epsilon$ and take the limit as $\epsilon \to 0^{+}$.
This provides a~discontinuity condition for the derivative of $v_l(\theta,\theta')$ with respect to~$\theta$ evaluated at
$\theta=\theta'$, namely
\begin{gather}
\label{eq:limv}
\lim\limits_{\epsilon \to 0^{+}} \left.
\frac{\partial v_l(\theta,\theta')}{\partial\theta} \right|_{\theta'-\epsilon}^{\theta'+\epsilon}=\frac{-1}{R^{d-2}}.
\end{gather}
After inserting~\eqref{eq:ul3} with~\eqref{eq:vl} into~\eqref{eq:limv}, substituting $z=\cos\theta'$, evaluating at $\theta=\theta'$,
and using~\eqref{WPQ},~\eqref{WPP}, we obtain $C=(-1)^{d/2-1+l}R^{2-d}$, $D=\frac12\pi(-1)^{(d-3)/2+l}R^{2-d}$.  Finally
\begin{gather*}
u_l(\theta,\theta')=
\begin{cases}
\displaystyle
\frac{(-1)^{d/2-1+l}}{R^{d-2}(\sin\theta\sin\theta')^{d/2-1}} \mathsf{P}_{d/2-1}^{-(d/2-1+l)} (\cos\theta_{<})
\mathsf{Q}_{d/2-1}^{d/2-1+l} (\cos\theta_{>})
 & \text{if}\  d \ \text{is even},
\vspace{1mm}\\
\displaystyle
\frac{\pi(-1)^{(d-3)/2+l}} {2R^{d-2}(\sin\theta\sin\theta')^{d/2-1}} \mathsf{P}_{d/2-1}^{-(d/2-1+l)} (\cos\theta_{<})
\mathsf{P}_{d/2-1}^{d/2-1+l} (\cos\theta_{>})
 & \text{if}\  d \ \text{is odd},
\end{cases}
\end{gather*}
and therefore through~\eqref{eq:SDY3}, this completes the proof.
\end{proof}

\subsection[Addition theorem for the azimuthal Fourier coef\/f\/icient on $\mathbf{S}_R^3$]{Addition theorem
for the azimuthal Fourier coef\/f\/icient on $\boldsymbol{\mathbf{S}_R^3}$}

As in~\cite[Section~5.1]{CohlKalII}, one may compute multi-summation addition theorems
for the azimuthal Fourier coef\/f\/icients on $\mathbf{S}_R^d$ by relating the azimuthal Fourier coef\/f\/icient
of $\mathcal{G}_R^d({\bf x},{{\bf x}^\prime})$ def\/ined
by~\eqref{FouriercoeffSRd} to its Gegenbauer polynomial expansion~\eqref{eq:SDY4}.
This is accomplished by using the addition theorem for hyperspherical harmonics~\eqref{eq:SUMY2} which expresses the Gegenbauer
polynomial with degree $l\in{\mathbf N}_0$ in~\eqref{eq:SDY4} as a~sum over the space of degenerate quantum numbers, of a~multiplicative
product of separated harmonics as a~function of geodesic polar coordinates on the hypersphere, such as~\eqref{eq:polar}.
For $d=3$ one obtains a~single-summation addition theorem using the notation~\eqref{alessgtr}.

\begin{corollary}
\label{thecorollary}
Let $m\in{\mathbf Z}$, $\theta,\theta',\theta_2,\theta_2'\in[0,\pi]$.
Then the azimuthal Fourier coefficient of a~fundamental solution of Laplace's equation on ${\mathbf S}_R^3$ is given~by
\begin{gather}
\mathsf{G}_m^{1/2}(\theta,\theta',\theta_2,\theta_2')= \frac{\pi\epsilon_m}{2\sqrt{\sin\theta\sin\theta'}}
\label{addthmfour}
\\
\qquad{}
\times\sum\limits_{l=|m|}^\infty (-1)^l(2l+l)\frac{(l-m)!}{(l+m)!} \mathsf{P}_l^m(\cos\theta_2)
\mathsf{P}_l^m(\cos\theta_2') \mathsf{P}_{1/2}^{-(1/2+l)}(\cos\theta_<) \mathsf{P}_{1/2}^{1/2+l}(\cos\theta_>).\nonumber
\end{gather}
\end{corollary}

\begin{proof}
By expressing~\eqref{eq:SDY4} for $d=3$, we obtain
\begin{gather}
\mathcal{G}_R^3({\bf x},{{\bf x}^\prime})=\frac{1}{8R\sqrt{\sin\theta\sin\theta'}} \sum\limits_{l=0}^\infty(-1)^l(2l\!+\!1)
\mathsf{P}_{1/2}^{-(1/2+l)}(\cos\theta_<) \mathsf{P}_{1/2}^{1/2+l}(\cos\theta_>) P_l(\cos\gamma),\!\!\!
\label{fundsolgegd3}
\end{gather}
since the Gegenbauer polynomial $C_l^{1/2}$ is the Legendre polynomial $P_l$.
By using the addition theorem for hyperspherical harmonics~\eqref{eq:SUMY2} with $d=3$ on has
\begin{gather}
P_l(\cos\gamma)=\sum\limits_{m=-l}^l\frac{(l-m)!}{(l+m)!} \mathsf{P}_{l}^{m}(\cos\theta_2) \mathsf{P}_{l}^{m}(\cos\theta_2')
e^{im(\phi-\phi')}.
\label{addthmsph}
\end{gather}
By inserting~\eqref{addthmsph} in~\eqref{fundsolgegd3} and reversing the order of the two summation symbols, we can compare the
resulting expression with~\eqref{Fourierexpansiond3} obtaining an expression for $\mathsf{G}_m^{1/2}:[0,\pi]^4\to{\mathbf R}$, in
standard hyperspherical coordinates~\eqref{eq:polar}.
This completes the proof of the corollary.
\end{proof}

The above addition theorem expresses the azimuthal Fourier coef\/f\/icients $\mathsf{G}_m^{1/2}$ given~by
Theo\-rem~\ref{azimuthalFourierexp3D} as the inf\/inite sum contribution from all meridional modes for a~given~$m$ value.
Alternatively, Theorem~\ref{azimuthalFourierexp3D} gives the azimuthal Fourier coef\/f\/icients as a~f\/inite sum of complete elliptic
integrals of the f\/irst, second and third kind which encapsulate the contribution from the inf\/inite meridional sum given~by
Corollary~\ref{thecorollary}.
By truncating the inf\/inite sum over the meridional modes, one obtains what is referred to in Euclidean space, as a~multipole method for
the potential problem.
The multipole method is well-known to be slowly convergent, even for moderate departure from spherical symmetry.
On the other hand, even for extremely deformed axisymmetric density distributions, Theorem~\ref{azimuthalFourierexp3D} gives the exact
solution to the potential problem in a~single~$m$ term! Note that the addition theorem~\eqref{addthmfour} reduces to the corresponding
result~\cite[equation~(2.4)]{CRTB} in the f\/lat-space limit.

\subsection{Spherically symmetric contribution to a~fundamental solution}

In the next section, we will take advantage of the following result to evaluate Newtonian potentials for spherically symmetric density
distributions.

\begin{corollary}
\label{sphsymfundsol}
Let $R>0$, $d=2,3,\ldots$, ${\mathcal H}_R^d:[0,\pi]\to{\mathbf R}$ be defined such that ${\mathcal H}_R^d(\Theta(\widehat{\bf
x},{\widehat{\bf x}^\prime})):= {\mathcal G}_R^d({\bf x},{{\bf x}^\prime})$ with geodesic distances measured from the origin.
Then
\begin{gather*}
{\mathcal G}_R^d({\bf x},{{\bf x}^\prime})\big|_{l=0} ={\mathcal H}_R^d(\theta_>).
\end{gather*}
Note that we are using the notation~\eqref{alessgtr}.
\end{corollary}

\begin{proof}
In dimension $d=2$ the result immediately follows from Theorem~\ref{azimuthalFourierexpansionS2}.
For $d\ge 3$ odd, using Theorem~\ref{Gegexpansion}, one has
\begin{gather*}
 {\mathcal G}_R^d({\bf x},{{\bf x}^\prime})\big|_{l=0} =\frac{(-1)^{(d-3)/2}\Gamma\left(\frac{d}{2}\right)}
{4\pi^{d/2-1}R^{d-2}(\sin\theta\sin\theta')^{d/2-1}} \mathsf{P}_{d/2-1}^{-(d/2-1)} (\cos\theta_{<}) \mathsf{P}_{d/2-1}^{d/2-1}
(\cos\theta_{>}),
\end{gather*}
and when combined with~\cite[equations~(14.5.18), (14.9.2)]{NIST}, the result follows. For $d\ge 4$ even, using Theorem~\ref{Gegexpansion}, one has
\begin{gather*}
 {\mathcal G}_R^d({\bf x},{{\bf x}^\prime})\big|_{l=0} =\frac{(-1)^{d/2-1}\Gamma\left(\frac{d}{2}\right)}
{2\pi^{d/2}R^{d-2}(\sin\theta\sin\theta')^{d/2-1}} \mathsf{P}_{d/2-1}^{-(d/2-1)} (\cos\theta_{<}) \mathsf{Q}_{d/2-1}^{d/2-1}
(\cos\theta_{>}),
\end{gather*}
which when combined with~\cite[equations~(14.5.18), (14.9.1)]{NIST} completes the proof. \end{proof}

\section{Applications}%\label{Applications}

The expansion formulas that we derive in the previous sections directly provide a~mechanism to obtain rapidly convergent solutions to
Newtonian potential problems in~$d$-dimensional hyperspherical geometry.
One may use these expansion formulae to obtain exact solutions either through numerical computation, or through symbolic derivation,
taking advantage of symmetries in density distributions.

An anonymous referee has suggested the following in regard to applications.
One potential source for applications of explicit fundamental solutions on hyperspheres is the mathematical geosciences (see for
instance~\cite{FreedenSchreiner}).
Furthermore, fundamental solutions of elliptic partial dif\/ferential equations on hyperspheres are often used as tools in numerical
analysis under the heading of spherical basis functions (SBFs) and that explicit formulas and series representations are essential for
practical implementations.
The use of SBFs in data-f\/itting and numerical solution of partial dif\/ferential equations (via collocation and Galerkin methods) is
plentiful (see for instance~\cite{Mhaskaretal2010,Narcowichetal2014}).

\subsection{Applications to Newtonian potential theory}\label{ApplicationsinPotentialTheory}

In this section we derive hyperspherical counterparts to some of the most elementary problems in classical Euclidean Newtonian potential
theory.

Consider Poisson's equation $-\Delta\Phi=\rho$ in Euclidean space~${\mathbf R}^d$ with Cartesian coordinates $(x_1,\ldots,x_d)$ and
$\Delta:=\sum\limits_{n=1}^d\frac{\partial^2}{\partial x_n^2}$.
One then has for $\rho:{\mathbf R}^d\to{\mathbf R}$ integrable (or even more generally $\rho\in({\mathcal D}({\mathbf R}^d))'$) and
$\Phi\in{C}^2({\mathbf R}^d)$, an integral solution to Poisson's equation
\begin{gather}
\label{NewtonRd}
\Phi({\bf x})=\int_{{\mathbf R}^d}{\mathcal N}^d({\bf x},{{\bf x}^\prime})\rho({\bf x}) d{{\bf x}^\prime},
\end{gather}
where ${\mathcal N}^d$ is given by~\eqref{eq:bolG}.
Furthermore, the total Newtonian binding energy $\mathcal E\in{\mathbf R}$, for a~density distribution~$\rho$ is given~by
(cf.~\cite[equation~(16), p.~64]{Chandrasekhar})
\begin{gather}
\mathcal E=\frac12\int_{{\mathbf R}^d}\rho({\bf x})\Phi({\bf x}) d{\bf x}.
\label{totalbinding}
\end{gather}
Now consider Poisson's equation on ${\mathbf S}_R^d$ which is given~by
\begin{gather}
-\Delta \Phi({\bf x})=\rho({\bf x}),
\label{Poissoneq}
\end{gather}
where $\rho:{\mathbf S}_R^d\to{\mathbf R}$ integrable (or even more generally $\rho\in({\mathcal D}({\mathbf S}_R^d))'$), $\Phi\in
C^2({\mathbf S}_R^d)$, and~$\Delta$ is given as in~\eqref{eq:delf}.
An integral solution to~\eqref{Poissoneq} is given through~\eqref{fundsoldefnJ},~\eqref{eq:vol} as
\begin{gather}
\Phi({\bf x})=\int_{{\mathbf S}_R^d}{\cal G}_R^d({\bf x},{{\bf x}^\prime})\rho({{\bf x}^\prime})d\vol'_g,
\label{Newton}
\end{gather}
and the total Newtonian binding energy~\eqref{totalbinding} for~$\rho$ is given analogously as
\begin{gather*}
\mathcal E =\frac12\int_{{\mathbf S}_R^d}\rho({\bf x})\Phi({\bf x}) d\vol_g.
%\label{binding}
\end{gather*}
Using our above derived azimuthal Fourier and Gegenbauer expansions for ${\cal G}_R^d$, we can exploit symmetries in~$\rho$ to derive
closed-form expressions for~$\Phi$ and $\mathcal E$ which reduce to their counterparts in Euclidean space.

The following three examples show direct applications for the expansions of a~fundamental solution of Laplace's equation
on hyperspheres given in this paper.
The f\/irst example shows an application of the azimuthal Fourier expansion for a~fundamental solution of Laplace's equation
on ${\mathbf S}_R^2$, namely Theorem~\ref{azimuthalFourierexpansionS2}.
The second example shows an application for the azimuthal Fourier expansion for a~fundamental solution of Laplace's equation
on ${\mathbf S}_R^3$, which is derived from Theorem~\ref{azimuthalFourierexp3D}.
The third example shows an application for the Gegenbauer expansion for a~fundamental solution of Laplace's equation on ${\mathbf S}_R^3$
which is derived from Theorem~\ref{Gegexpansion}.

\begin{example}
\label{2ballexample}
Newtonian potential of a~uniform density 2-ball with geodesic radius $R\theta_0$ in ${\mathbf S}_R^2$ (hereafter, uniform density
2-disc).
Consider a~uniform density 2-disc in ${\mathbf S}_R^2$ def\/ined using~\eqref{eq:polar} for all $\phi\in[-\pi,\pi)$ such that
\begin{gather}
\label{uniformballdensity23}
\rho(\mathbf{x}):= \begin{cases}
\rho_0 & \text{if}\quad \theta\in[0,\theta_0],
\\
0 & \text{if}\quad \theta\in(\theta_0,\pi],
\end{cases}
\end{gather}
where $|\rho_0|>0$.
Due to the axial symmetry in~$\rho$ (invariance under rotations centered about the origin of ${\mathbf S}_R^2$ at $\theta=0$), the only
non-zero contribution to~$\Phi$ is from the $n=0$ term in Theorem~\ref{azimuthalFourierexpansionS2}.
This term is given by Corollary~\ref{sphsymfundsol} and~\eqref{fs2} as
\begin{gather}
{\cal G}_R^2({\bf x},{{\bf x}^\prime})\big|_{n=0}=\frac{1}{2\pi}\log\cot \frac{\theta_>}{2}.
\label{Greenn0}
\end{gather}
Applying~\eqref{Greenn0} in~\eqref{Newton} with~\eqref{uniformballdensity23}, using elementary trigonometric integration, one obtains
\begin{gather}
\label{uniformballpotential2}
\Phi(\mathbf{x}):=
\begin{cases}
\displaystyle \rho_0R^2\left(\log\cot\frac{\theta}{2}-\cos\theta_0\log\cot\frac{\theta_0}{2}+\log\frac{\sin\theta}{\sin\theta_0}\right)
 & \text{if}\quad \theta\in[0,\theta_0],
\\
\displaystyle \rho_0R^2(1-\cos\theta_0)\log\cot\frac{\theta}{2} & \text{if}\quad \theta\in(\theta_0,\pi],
\end{cases}
\end{gather}
for all $\phi\in[-\pi,\pi)$.

In comparison, consider the problem to obtain the Newtonian potential for an uniform density 2-disc with radius $r_0$ in ${\mathbf
R}^2$, with points parametrized using polar coordinates $(r\cos\phi,r\sin\phi)$.
Using~\eqref{2discpotential} we see that ${\mathcal N}^2({\bf x},{{\bf x}^\prime})\bigr|_{n=0}=-\frac{1}{2\pi}\log r_>$, which yields
through~\eqref{NewtonRd}
\begin{gather}
\label{uniformballpotentialR2}
\Phi(\mathbf{x}):= \begin{cases}
-\dfrac{\rho_0}{4}\big(r^2-r_0^2+2r_0^2\log r_0\big) & \text{if}\quad r\in[0,r_0],
\\
-\dfrac12\rho_0r_0^2\log r & \text{if}\quad r\in(r_0,\infty),
\end{cases}
\end{gather}
for all $\phi\in[-\pi,\pi)$.
The potential on the uniform density 2-disc~\eqref{uniformballpotential2} reduces to~\eqref{uniformballpotentialR2} in the Euclidean space.
This is accomplished by performing an asymptotic expansion as $\theta,\theta_0\to 0^{+}$ in the f\/lat-space limit $R\to\infty$.
Note however that in this two-dimensional example, such as in other two-dimensional problems (see for example~\cite[p.~151]{Fol4}),
we must subtract a~constant ($c_R:=\frac12\rho_0r_0^2\log(2R)$) which tends to inf\/inity.
This is required in the f\/lat-space limit $R\to\infty$ for~$\Phi$ in~\eqref{uniformballpotential2} to obtain the f\/inite
limit~\eqref{uniformballpotentialR2}.
The potential~$\Phi$ satisf\/ies~\eqref{Poissoneq}, and therefore, so does any solution $\Phi+c_R$, with $c_R$ constant.
Using the geodesic distances $r\sim R\theta$, $r_0\sim R\theta_0$, the result follows using Maclaurin expansions for the logarithmic and
trigonometric functions in~\eqref{uniformballpotential2}.

For the uniform density 2-disc, the total Newtonian binding energy is
\begin{gather}
\mathcal E_{2\text{-disc}}=\frac{\pi}{4}\rho_0^2R^4\bigl((1-4\cos\theta_0+\cos(2\theta_0))\log\cot\frac{\theta_0}{2}-4\log\cos\frac{\theta_0}{2}
\nonumber
\\
\phantom{\mathcal E_{2\text{-disc}}=}
-2\log\sin\theta_0 +2\cos(\theta_0) +2(\log 2-1) \bigr).
\label{binding2}
\end{gather}
For the uniform density 2-disc embedded in ${\mathbf R}^2$, the total Newtonian binding energy is
\begin{gather}
\mathcal E_{{\mathbf R}^2\text{-}2\text{-disc}}=\frac{\pi}{16}\rho_0^2r_0^4 \left[-4\log r_0+1
+\frac{r_0^2}{6R^2}\left(2\log\frac{r_0}{2}-1\right) \right],
\label{EuclideanE2}
\end{gather}
where we have included the lowest-order correction term due to curvature originating from~\eqref{binding2}.
In the f\/lat-space limit $R\to\infty$,~\eqref{binding2} reduces to~\eqref{EuclideanE2}, with the term corresponding to $c_R$ being
subtracted of\/f.
\end{example}

\begin{example}
Newtonian potential of a~uniform density circular curve segment with geodesic length $2R\varphi$ with $\varphi\in(0,\pi]$ in ${\mathbf
S}_R^3$.
Consider a~uniform density circular curve segment using Hopf coordinates~\eqref{eq:hopf} on ${\mathbf S}_R^3$.
Recall $\vartheta\in[0,\frac{\pi}{2}]$, $\phi_1,\phi_2\in[-\pi,\pi)$.
The density distribution~$\rho$ is given by the uniform density circular curve segment
\begin{gather}
\label{uniformcurve}
\rho(\mathbf{x}):=
\begin{cases}
\dfrac{\rho_0\delta(\vartheta)}{R^2\sin\vartheta\cos\vartheta} & \text{if}\quad \phi_1\in[-\varphi,\varphi),
\\
0 & \text{otherwise}.
\end{cases}
\end{gather}
This uniform density distribution is rotationally-invariant about the $\phi_2$ angle.
In order to make a~connection with the corresponding Euclidean space problem, consider the f\/lat-space limit $R\to\infty$ of the density
distribution.
Recall from Section~\ref{Theflatspacelimit} that $\vartheta\sim \frac{r}{R}$, $\phi_1\sim\frac{z}{R}$, with $\varphi\sim\frac{L}{R}$,
$L>0$.
One has since $\delta(ax)=\frac{\delta(x)}{|a|}$, that
\begin{gather*}
\rho({\bf x})=\frac{\rho_0\delta(\vartheta)}{R^2\sin\vartheta\cos\vartheta} \sim
\frac{\rho_0\delta(\frac{r}{R})}{R^2\sin\left(\frac{r}{R}\right) \cos\left(\frac{r}{R}\right)} \sim \frac{\rho_0\delta(\frac{r}{R})}{Rr}
= \frac{\rho_0\delta(r)}{r},
\end{gather*}
for all $z\in[-L,L]$.
This is the axisymmetric density distribution of a~uniform density line segment centered on the origin, along the~$z$-axis in ${\mathbf R}^3$.

In order to solve for the Newtonian potential of this density distribution~\eqref{uniformcurve} on ${\mathbf S}_R^3$ using Hopf
coordinates, one must obtain the Riemannian structure, namely
\begin{gather*}
ds^2=R^2\big(d\vartheta^2+\cos^2\vartheta d\phi_1^2+\sin^2\vartheta d\phi_2^2\big).
\end{gather*}
This yields the metric coef\/f\/icients and therefore in Hopf coordinates on ${\mathbf S}_R^3$ one has the Riemannian volume
measure~\cite[p.~29]{Lee}
\begin{gather*}
d\vol_g=R^3\cos\vartheta\sin\vartheta d\vartheta d\phi_1 d\phi_2.
\end{gather*}
and the geodesic distance through~\eqref{eq:cosd}
\begin{gather*}
d({\bf x},{{\bf x}^\prime})=R\cos^{-1}\left(\cos\vartheta\cos\vartheta'\cos(\phi_1-\phi_1')
+\sin\vartheta\sin\vartheta'\cos(\phi_2-\phi_2')\right).
\end{gather*}
A~fundamental solution of Laplace's equation is given through~\eqref{Fourierexpansiond3} as
\begin{gather*}
{\cal G}_R^3({\bf x},{{\bf x}^\prime})=\frac{A/B+\cos\psi} {4\pi R \sqrt{\left(\frac{1-A}{B}-\cos\psi\right)
\left(\frac{1+A}{B}+\cos\psi\right)}},
\end{gather*}
where $A:[0,\frac{\pi}{2}]^2\times[-\pi,\pi)^2\to[0,1]$, $B:[0,\frac{\pi}{2}]^2\to[0,1]$, are def\/ined as
\begin{gather*}
A(\vartheta,\vartheta',\phi_1,\phi_1'):=\cos\vartheta\cos\vartheta'\cos(\phi_1-\phi_1'),
\qquad
B(\vartheta,\vartheta'):=\sin\vartheta\sin\vartheta'.
\end{gather*}
The azimuthal Fourier coef\/f\/icients of $\mathfrak{g}^3=4\pi R{\cal G}_R^3$ are therefore given as an integral of the same form as that
given in standard hyperspherical coordinates, cf.~\eqref{eq:Smhalf1},
\begin{gather*}
\mathsf{G}_m^{1/2}(\vartheta,\vartheta',\phi_1,\phi_1')= \frac{\epsilon_m}{\pi} \int_{-1}^{1} \frac{\left(x+A/B \right)
T_m(x) dx}{\sqrt{(1-x^2)\left(\frac{1-A}{B}-x\right) \left(\frac{1+A}{B}+x\right)}},
\end{gather*}
where $x=\cos(\phi_2-\phi_2')$.
Since our density distribution is axisymmetric about the $\phi_2$ angle, the only contribution to the Newtonian potential originates
from the $m=0$ term.
We can evaluate this term using the methods of Section~\ref{four3} which results in
\begin{gather*}
\mathsf{G}_0^{1/2}(\vartheta,\vartheta',\phi_1,\phi_1')=
\frac{2}{\pi}\left[\frac{K(k)}{\sqrt{(1-A+B)(1+A+B)}}+\frac{(A+B-1)\Pi(\alpha^2,k)}{\sqrt{(A+B+1)(1-A+B)}} \right].
\end{gather*}

The Newtonian potential for the density distribution~\eqref{uniformcurve} is given through~\eqref{Newton}.
The $\vartheta'$ integration selects $\vartheta'=0$ in the integrand which implies that $A=\cos\vartheta\cos(\phi_1-\phi_1')$, $B=0$,
$k=\alpha^2=0$, and since $K(0)=\Pi(0,0)=\frac{\pi}{2}$,
\begin{gather*}
\mathsf{G}_0^{1/2}(\vartheta,0,\phi_1,\phi_1')= \frac{A}{\sqrt{1-A^2}}=\frac{\cot\vartheta\cos(\phi_1-\phi_1')}
{\sqrt{1+\cot^2\vartheta\sin^2(\phi_1-\phi_1')}}.
\end{gather*}
The Newtonian potential of the uniform density circular curve segment is given~by
\begin{gather*}
\Phi({\bf x})=\frac{\rho_0\cot\vartheta}{2} \int_{-\varphi}^\varphi \frac{\cos(\phi_1-\phi_1')d\phi_1'}
{\sqrt{1+\cot^2\vartheta\sin^2(\phi_1-\phi_1')}}.
\end{gather*}
With the substitution $u=\cot\vartheta\sin(\phi_1-\phi_1')$, one f\/inally has
\begin{gather*}
\Phi({\bf x})=\frac{\rho_0}{2} \int_{\cot\vartheta\sin(\phi_1-\varphi)}^{\cot\vartheta\sin(\phi_1+\varphi)}
\frac{du}{\sqrt{1+u^2}}=\frac{\rho_0}{2}\log\left| \frac {\sin(\phi_1+\varphi)+\sqrt{\tan^2\vartheta+\sin^2(\phi_1+\varphi)}}
{\sin(\phi_1-\varphi)+\sqrt{\tan^2\vartheta-\sin^2(\phi_1+\varphi)}} \right|,
\end{gather*}
since the antiderivative of $(1+z^2)^{-1/2}$ is $\sinh^{-1}z=\log|z+\sqrt{z^2+1}|$. In the f\/lat-space limit $R\to\infty$ we have
\begin{gather*}
\Phi({\bf x})=\frac{\rho_0}{2} \log\left| \frac {z+L+\sqrt{r^2+(z+L)^2}} {z-L+\sqrt{r^2+(z-L)^2}} \right|,
\end{gather*}
which is the potential in cylindrical coordinates of a~uniform density line segment on the $z'$-axis with $z'\in[-L,L]$.
\end{example}

\begin{example}
Newtonian potential of a~uniform density 3-ball with geodesic radius $R\theta_0$ in~${\mathbf S}_R^3$ (hereafter, uniform density 3-ball).
Consider a~uniform density 3-ball def\/ined using~\eqref{eq:polar} for all $\theta_2\in[0,\pi]$, $\phi\in[-\pi,\pi)$ such
that~\eqref{uniformballdensity23} is true with $|\rho_0|>0$.
Due to the spherical symmetry in~$\rho$ (invariance under rotations centered about the origin of ${\mathbf S}_R^3$ at $\theta=0$), the
only non-zero contribution to~$\Phi$ is from the $l=0$ term in Theorem~\ref{Gegexpansion}.
This term is given through Corollary~\ref{sphsymfundsol} and~\eqref{Fourierexpansiond3} as
\begin{gather}
{\cal G}_R^3({\bf x},{{\bf x}^\prime})\bigr|_{l=0} =\frac{1}{4\pi R}\cot\theta_>.
\label{Greenl0}
\end{gather}
Applying~\eqref{Greenl0} in~\eqref{Newton} with~\eqref{uniformballdensity23}, using elementary trigonometric integration, one obtains
\begin{gather}
\label{uniformballpotential3}
\Phi(\mathbf{x}):=
\begin{cases}
\dfrac{\rho_0R^2}{2}[\theta\cot\theta-\cos^2\theta_0] & \text{if}\quad \theta\in[0,\theta_0],
\vspace{1mm}\\
\dfrac{\rho_0R^2}{2}[\theta_0-\sin\theta_0\cos\theta_0]\cot\theta & \text{if}\quad \theta\in(\theta_0,\pi],
\end{cases}
\end{gather}
for all $\theta_2\in[0,\pi]$, $\phi\in[-\pi,\pi)$.

In comparison, consider the well-known freshman problem to obtain the Newtonian potential for a~uniform density 3-ball with radius $r_0$
in ${\mathbf R}^3$.
Using the generating function for Gegenbauer polynomials~\cite[equation~(18.12.4)]{NIST} for $d=3$
one obtains ${\mathcal N}^3({\bf x},{{\bf x}^\prime})\bigr|_{l=0}=\frac{1}{4\pi r_>}$, which yields
\begin{gather}
\label{uniformballpotentialR3}
\Phi(\mathbf{x}):= \begin{cases}
\dfrac{\rho_0}{6}\big(3r_0^2-r^2\big) & \text{if}\quad r\in[0,r_0],
\vspace{1mm}\\
\dfrac{\rho_0r_0^3}{3r} & \text{if}\quad r\in(r_0,\infty).
\end{cases}
\end{gather}
The Newtonian potential of the 3-ball~\eqref{uniformballpotential3} reduces to~\eqref{uniformballpotentialR3} in the f\/lat-space limit,
namely $\theta,\theta_0\to 0^{+}$ in the f\/lat-space limit $R\to\infty$.
Using the geodesic distances $r\sim R\theta$, $r_0\sim R\theta_0$, the result follows using the f\/irst two terms of Laurent series for
$\cot\theta$.

For the uniform density 3-ball, the total Newtonian binding energy is
\begin{gather}
\mathcal E_{3\text{-ball}}=\frac{\pi}{16}\rho_0^2R^5
\left(-4\theta_0+4\sin(2\theta_0)+\sin(4\theta_0)-8\theta_0\cos(2\theta_0)\right).
\label{Esphere3}
\end{gather}
Furthermore, the classic result for the 3-ball embedded in ${\mathbf R}^3$ is given through~\eqref{uniformballpotentialR3}, namely
\begin{gather}
\mathcal E_{{\mathbf R}^3\text{-}3\text{-ball}} =\frac{4\pi}{15}\rho_0^2r_0^5\left(1-\frac{13}{21}\frac{r_0^2}{R^2}\right)
=\frac{3M^2}{20\pi r_0}\left(1-\frac{13}{21}\frac{r_0^2}{R^2}\right),
\label{EuclideanE3}
\end{gather}
where we have included the lowest-order correction term due to curvature. Since the total mass of the 3-ball
is given by $M=\frac{4}{3}\pi \rho_0 r_0^3$, the standard result $(-3GM^2/(5r_0))$ from Newtonian
gravity is obtained by mapping $M\mapsto-4\pi GM$ for the~$M$ which originated from the Newtonian potential (by solving the Poisson
equation $\Delta\Phi=4\pi G\rho$).
If one performs a~Maclaurin expansion about $\theta_0$ in~\eqref{Esphere3}, then we see that
\begin{gather*}
-4\theta_0+4\sin(2\theta_0)+\sin(4\theta_0)-8\theta_0\cos(2\theta_0)= \frac{64}{15}\theta_0^5-\frac{832}{315}\theta_0^7+O\big(\theta_0^9\big),
\end{gather*}
and since $\theta_0\sim r_0/R$, then~\eqref{EuclideanE3} follows in the f\/lat-space limit $R\to\infty$.
\end{example}

\subsection{Applications to superintegrable potentials}\label{Applicationstosuperintegrable}

Superintegrable potentials, and in fact any potentials such that, $\Phi\in{C}^2({\mathbf S}_R^d)$ will satisfy Poisson's equation
$-\Delta \Phi({\bf x})=\rho({\bf x})$, for a~particular paired density distribution $\rho:{\mathbf S}_R^d\to{\mathbf R}$.
Therefore, given a~superintegrable $\Phi({\bf x})$ on ${\mathbf S}_R^d$, one may obtain a~density distribution corresponding to the
superintegrable potential by evaluating the negative Laplace--Beltrami operator on the~$d$-dimensional hypersphere acting on the
superintegrable potential.
Furthermore, one may then re-obtain the superintegrable potential from $\rho({\bf x})$ using~\eqref{Newton}.
Fourier and Gegenbauer series expansions of fundamental solutions, give you the ability to demonstrate this connection when
a~superintegrable potential and the corresponding density distribution satisfy a~degree of symmetry, such as for rotationally-invariant
or spherically-symmetric density distributions.

Below we present two isotropic examples of this correspondence on ${\mathbf S}_R^d$.
However, such paired density distributions are available for all superintegrable potentials, no matter how complicated, by simply
applying the negative Laplace--Beltrami operator to those potentials.
Therefore a~superintegrable potential can be re-obtained by convolution with the paired density distribution and a~fundamental solution
for Laplace's equation on the~$d$-dimensional~$R$-radius hypersphere.

\begin{example}
Take for instance the superintegrable isotropic oscillator potential on ${\mathbf S}_R^d$ in standard hyperspherical
coordinates~\eqref{eq:polar}, given by~\cite[equation~(1.4)]{HerranzBallesteros}
\begin{gather}
\Phi({\bf x})=\alpha\tan^2\theta,
\label{isooscillator}
\end{gather}
for $|\alpha|>0$.
In this case, one obtains through~\eqref{eq:delf}, the paired density distribution
\begin{gather*}
\rho({\bf x})=\frac{-2\alpha}{R^2} \big(1+\tan^2\theta\big)\big(3\tan^2\theta+d\big).
\end{gather*}
Note that this density distribution has an inf\/inite total mass given~by
\begin{gather*}
M=\int_{{\mathbf S}_R^d}\rho({\bf x})d\vol_g,
\end{gather*}
since through~\eqref{eq:vol} one may obtain
\begin{gather}
\label{hypersurfarea}
\int_{{\mathbf S}_R^{d-1}}d\omega=\frac{2\pi^{d/2}}{\Gamma\left(\frac{d}{2}\right)},
\end{gather}
and also using that $\int_0^\pi(1+\tan^2\theta)(3\tan^2\theta+d)\sin^{d-1}\theta d\theta$ diverges.

On the other hand, one may use~\eqref{Newton} make a~concrete correspondence between the superintegrable potential~\eqref{isooscillator}
and the convolution of a~fundamental solution of Laplace's equation in hyperspherical geometry.
Using~\eqref{Newton},~\eqref{hypersurfarea}, and Corollary~\ref{sphsymfundsol}, one has
\begin{gather*}
\tan^2\theta=\frac{-4\pi^{d/2} R^{d-2}}{\Gamma\left(\frac{d}2\right)} \int_{{\mathbf S}_R^d}{\mathcal H}_R^d(\theta_>)
\big(1+\tan^2\theta'\big)\big(3\tan^2\theta'+d\big) \sin^{d-1}\theta'd\theta',
\end{gather*}
which is equivalent to
\begin{gather*}
\tan^2\theta=\frac{-(d-2)!}{2^{d/2-2}\Gamma\left(\frac{d}2\right)} \int_0^\pi\frac{\sin^{d-1}\theta'}{\sin^{d/2-1}\theta_>}
\mathsf{Q}_{d/2-1}^{1-d/2}(\cos\theta_>) \big(1+\tan^2\theta'\big)\big(3\tan^2\theta'+d\big)d\theta'.
\end{gather*}
By evaluating the above integral, one arrives at a~new def\/inite integral valid for $d=2,3,\ldots$, namely
\begin{gather*}
\int_\theta^\pi \mathsf{Q}_{d/2-1}^{1-d/2}(\cos\theta') \big(1+\tan^2\theta'\big) \big(3\tan^2\theta'+d\big)\sin^{d/2}\theta'd\theta'
\\
\qquad{}
= -\tan\theta\big(1+\tan^2\theta\big)\sin^{d/2}\theta\mathsf{Q}_{d/2-1}^{1-d/2}(\cos\theta)-\frac{2^{d/2-2}\Gamma\left(\frac{d}2\right)\tan^2\theta}{(d-2)!}.
\end{gather*}

\end{example}

\begin{example}
Take for instance the superintegrable Kepler--Coulomb potential on ${\mathbf S}_R^d$ in standard hyperspherical
coordinates~\eqref{eq:polar}, also given by~\cite[equation~(1.4)]{HerranzBallesteros}
\begin{gather}
\Phi({\bf x})=-\alpha\cot\theta,
\label{isokepler}
\end{gather}
for $|\alpha|>0$.
In this case, one obtains through~\eqref{eq:delf}, the paired density distribution
\begin{gather*}
\rho({\bf x})=\frac{(3-d)\alpha}{R^2}\cot\theta\big(1+\cot^2\theta\big).
\end{gather*}
Note that this density distribution has an zero total mass for all dimensions greater than two, except for $d=3$, which gives a~f\/inite
total mass.
To prove that the total mass vanishes for $d=2$, one may approach the endpoints of integration in a~limiting fashion.
For $d=3$, this superintegrable potential is a~fundamental solution for Laplace's equation on the 3-sphere~\eqref{Fourierexpansiond3}
(see for instance~\cite{Schrodinger40}), and therefore the paired density distribution is proportional to the Dirac delta distribution
which integrates to unity.

Use~\eqref{Newton} to make a~concrete correspondence between the superintegrable potential~\eqref{isokepler} and the convolution of
a~fundamental solution of Laplace's equation in hyperspherical geometry.
Using~\eqref{Newton},~\eqref{hypersurfarea}, and Corollary~\ref{sphsymfundsol}, one has
\begin{gather*}
-\cot\theta=\frac{(3-d)(d-2)!}{2^{d/2-1}\Gamma\left(\frac{d}2\right)} \int_0^\pi \frac{\sin^{d-1}\theta'}{\sin^{d/2-1}\theta_>}
\mathsf{Q}_{d/2-1}^{1-d/2}(\cos\theta_>) \big(1+\cot^2\theta'\big)\cot\theta'd\theta'.
\end{gather*}
By evaluating the above integral, one arrives at a~new def\/inite integral valid for $d=4,5,\ldots$, namely
\begin{gather*}
\int_\theta^\pi \mathsf{Q}_{d/2-1}^{1-d/2}(\cos\theta') \big(1+\cot^2\theta'\big) \cot\theta'\sin^{d/2}\theta'd\theta' \\
\qquad
= -\frac{\sin^{d/2-2}\theta}{(d-3)} \mathsf{Q}_{d/2-1}^{1-d/2}(\cos\theta) +\frac{2^{d/2-1}\Gamma\left(\frac{d}2\right)\cot\theta}{(d-3)(d-2)!}.
\end{gather*}
For $d=2$, in order to obtain the paired density distribution for~\eqref{isokepler}, one must approach the endpoints of integration in
a~limiting fashion.
By taking the integration over $\theta'\in[\epsilon,\pi-\epsilon]$, and then taking the limit $\epsilon\to 0^{+}$, the singular terms
which arise cancel.
Furthermore, one must add a~constant (whose gradient vanishes for all values of~$\epsilon$) which tends to inf\/inity as $\epsilon\to
0^{+}$ to~\eqref{isokepler}.
This is often the case in two-dimensions, see for instance Example~\ref{2ballexample}.
In this case, it is not dif\/f\/icult to see that one must also add a~Dirac delta distribution source term at the origin to the density
distribution.
The superintegrable density-potential pair for~\eqref{isokepler} on ${\mathbf S}_R^2$ is therefore
\begin{gather*}
\Phi({\bf x})=\alpha\cot\theta+\alpha\cot\epsilon+\frac{\alpha}{\sin\epsilon} \log\tan\frac{\epsilon}{2},
\end{gather*}
and
\begin{gather*}
\rho({\bf x})=-\frac{\alpha}{R^2}\cot\theta\big(1+\cot^2\theta\big)+\frac{2\alpha\delta(\theta)} {\sin\epsilon\sin\theta},
\end{gather*}
with~\eqref{Newton} being satisf\/ied as $\epsilon\to0^{+}$.
For $d=3$, no new information is gained since in this case the superintegrable potential~\eqref{isokepler} is simply a~fundamental
solution of Laplace's equation on~${\mathbf S}_R^3$~\eqref{Fourierexpansiond3}.
\end{example}

\subsection*{Acknowledgements}

Much thanks to Willard Miller and George Pogosyan for valuable
discussions. We would also like to express our gratitude to the anonymous referees and the editors for
this special issue in honour of Luc Vinet, for their signif\/icant contributions.

\pdfbookmark[1]{References}{ref}
\LastPageEnding

\end{document}